\documentclass[11pt,a4paper]{article}
\usepackage{amsmath, amsfonts}
  \usepackage{paralist}
  \usepackage{microtype}

\newtheorem{thm}{Theorem}[section]

\newtheorem{lem}[thm]{Lemma}

\newtheorem{definition}[thm]{Definition}

\numberwithin{equation}{section}

\begin{document}

\title{\bf Stability of normalized solitary waves for three coupled nonlinear Schr\"{o}dinger equations}
\author{\small SANTOSH BHATTARAI}

\date{}

\maketitle

\begin{footnotesize}
\noindent {\bf Abstract.} In this paper we establish existence and stability results concerning fully
nontrivial solitary-wave solutions to 3-coupled nonlinear Schr\"{o}dinger system
\begin{equation*}
 i\partial_t u_{j}+%
\partial_{xx}u_{j}+ \left(\sum_{k=1}^{3} a_{kj} |u_k|^{p}\right)|u_j|^{p-2}u_j = 0, \ j=1,2,3,
\end{equation*}
where $u_j$ are complex-valued functions of $(x,t)\in \mathbb{R}^{2}$ and $a_{kj}$ are
positive constants satisfying $a_{kj}=a_{jk}$ (symmetric attractive case).
Our approach improves many of the previous known results. In all methods used previously
to study solitary waves,
which we are aware of, the variational problem has consisted of
finding the extremum of an energy functional
subject to the constraints that were not independently chosen.
Here we study a problem of minimizing the energy functional subject to
three independent $L^2$ mass constraints and
establish existence and stability results for
a true three-parameter family of solitary waves.

\medskip

\end{footnotesize}

%The title of your section 1
\section{Introduction}

\noindent In recent years much attention has been given to the study of coupled
nonlinear Schr\"{o}dinger (CNLS) equations because of their
applications in a variety of physical and biological settings.
The CNLS equation models physical systems in which the field has more than one component.
For example, the CNLS equations play an important role in
wavelength-division multiplexing \cite{[Cha], [Mol]} and
multichannel bit-parallel- wavelength optical fiber networks \cite{[Yeh]}, where the pulses
propagate at least in two channels simultaneously. In addition, the CNLS equations arise in plasma physics \cite{[Som]},
multispecies and spinor Bose-Einstein condensates \cite{[Dal], [Ho], [Kawa],[Ohmi]},
biophysics \cite{[Sco]}, nonlinear Rossby waves \cite{[Sun]}, to name a few.

\medskip

\noindent In this paper, we consider the time-dependent 3-coupled nonlinear Schr\"{o}dinger equations given by
\begin{equation}\label{CNLS}
\left\{
\begin{aligned}
& iu_{j,t}+%
u_{j,xx}+ \left(\sum_{k=1}^{3} a_{kj} |u_k|^{p}\right)|u_j|^{p-2}u_j = 0, \ j=1,2,3, \\
& \ u_{j}=u_{j}(x,t)\in \mathbb{C},\ (x,t)\in \mathbb{R}^{2},
\end{aligned}
\right.
\end{equation}
where $u_j$ are dimensionless complex amplitude of the $j$-th component of the underlying physical system
 and $a_{jk}$ are constants satisfying $a_{jk}=a_{kj}$ for all $j,k\in \{1,2,3\}.$
The interaction matrix $(a_{jk})_{j,k=1}^{3}$ contains information about the nature of
the interactions between the different components of the wave functions. The Landau constants
$a_{jj}$ describe the self-modulation of the wave packets, and the coupling constants $a_{kj} \ (k\neq j)$ are the
wave-wave interaction coefficients, which describe the cross-modulation of the wave packets.
The interaction is (purely) attractive if all couplings are positive and the
interaction is (purely) repulsive when they all are negative.
Here we study the symmetric attractive interactions.
Throughout this paper, we shall denote by $u_{j, t}$ the partial derivative of $u_j$ with respect
to $t$ and by $u_{j, xx}$ the second partial derivative with respect to $x.$

\medskip

\noindent By a solitary-wave solution of \eqref{CNLS} we mean
a function $(\Phi_1,\Phi_2,\Phi_3)$ such that $\Phi_j$ are in $H^1(\mathbb{R})$ and $(u_1,u_2,u_3)$ defined by
\begin{equation} \label{SO}
u_{j}(x,t)=e^{i(\omega _{j}-\sigma ^{2})t+i\sigma x+i\beta_j}\Phi_{j} (x-2\sigma t),\ j=1,2,3,
\end{equation}
is a solution of \eqref{CNLS} for some real numbers
$\omega_{j}, \sigma,$ and $\beta_j.$
When $\sigma=0,$ solutions of the form \eqref{SO} are usually referred to
as standing-wave solutions.
Inserting solitary waves ansatz \eqref{SO} into \eqref{CNLS}, we
see that $(\Phi_{1}, \Phi_{2}, \Phi_{3})$ solves the following system of ordinary differential equations
\begin{equation}\label{ODE}
-\Phi _{j}^{\prime \prime }+\omega _{j}\Phi
_{j}= \left(\sum_{k=1}^{3}a_{jk}|\Phi _{k}|^{p}\right)|\Phi_{j}|^{p-2}\Phi_{j}, \ j=1,2,3.
\end{equation}
System \eqref{ODE} has many semi-trivial
(or collapsing) solutions, i.e., solution $(u_1,u_2,u_3)$ with at least one, but not all, component being zero.
In these cases the system collapses into system
with fewer components. A natural question relevant for 3-coupled nonlinear systems such as \eqref{ODE}
is the existence and stability results
of nontrivial solutions
(we will call a solution nontrivial if all three components
of the solution are non-zero). In the literature these solutions are also referred to as \textit{co-existing solutions}.
This paper aims to address the issues of existence of nontrivial solutions to \eqref{ODE} and
the stability of corresponding solitary waves
for the full equations \eqref{CNLS}.

\medskip

\noindent{In what follows we denote by $Y$ the product space $H^{1}(\mathbb{R})\times H^1(\mathbb{R})\times H^1(\mathbb{R}).$
The following definition of stability is used throughout the paper.}
\begin{definition}\label{def1}
Let
$\Sigma \subset Y$ be a set of vectors of solitary-wave profiles $\mathbf{\Phi}=(\Phi_1,\Phi_2,\Phi_3);$
i.e., each $\mathbf{\Phi} \in \Sigma$ corresponds to a
solution $\mathbf{u}(x,t)$ of \eqref{CNLS}.
We say that $\Sigma$ is a stable set of solitary-wave profiles if
for any $\epsilon>0$ there exists $\delta>0$ such that for every $\Psi$ (in a suitable space $X$ of initial data)
satisfying $\inf_{\mathbf{w} \in \Sigma}\|\Psi-\mathbf{w}\|_Y <\delta,$
the solution $\mathbf{u}(x,t)$ of \eqref{CNLS} with $\mathbf{u}(x,0)=\Psi(x)$ satisfies
 \begin{equation*}
\sup_{t\in \mathbb{R}}\inf_{\mathbf{w} \in \Sigma} \|\mathbf{u}(x,t)-\mathbf{w}\|_Y <\epsilon.
\end{equation*}
\end{definition}

\noindent Implicit in the above definition of stability is the assumption
that the initial-value problem associated to \eqref{CNLS}
is globally well-posed in some space $X$ of ordered triples of functions of $x.$
Here we adapt the standard notion of the well-posedness. More precisely, we say that the IVP for \eqref{CNLS} is
globally
well-posed (g.w.p.) in $X$ if for a given $\Psi \in X$ there exists a
unique $\mathbf{u}(x,t)$ such that $\mathbf{u}(x,0)=\Psi(x), \mathbf{u}(\cdot, t) \in X$ for all $t \in \mathbb{R},$
and $\mathbf{u}(x,t)$ solves \eqref{CNLS} in some (possibly weak) sense. Moreover, the map
$t\mapsto \mathbf{u}(\cdot,t)$ is in the space $\mathcal{C}(\mathbb{R}; X)$ of continuous maps from
$\mathbb{R}$ to $X,$ and the solution map $\Psi \mapsto \mathbf{u}(x,t)$ from the initial data
to the solution
defines a continuous map from
$X$ to $\mathcal{C}(\mathbb{R}; X)$ in the appropriate topology.
For our purposes, the well-posedness result in \cite{[Ngu3]} (see also \cite{[Ca]}) is most convenient
because it is set in the energy space $Y$ and their method works for the range $2\leq p <3.$
It has been proved in \cite{[Ngu3]} that for any initial data $\mathbf{u}(x,0)$ lying
in the
space $Y,$ there exists a unique solution $\mathbf{u}(x,t)$ of \eqref{CNLS} in $\mathcal{C}(\mathbb{R}, Y)$
emanating from $\mathbf{u}(x,0),$ and $\mathbf{u}(x,t)$ satisfies
\begin{equation*}
\mathcal{H}(\mathbf{u}(x,t))=\mathcal{H}(\mathbf{u}(x,0))\ \textrm{and}\ \mathcal{Q}(u_{j}(x,t))=\mathcal{Q}(u_{j}(x,0)),
\end{equation*}
where $\mathcal{H}$ and $\mathcal{Q}$ are the following conserved quantities
\begin{equation}\label{Edef}
\mathcal{H}(\mathbf{u}(x,t))=\int_{-\infty }^{\infty }\left(
\sum_{j=1}^{3}|%
u_{j, x}(x,t)|^{2}-\frac{1}{p}\sum_{k,j=1}^{3}a_{kj}
|u_{k}(x,t)|^{p}|u_{j}(x,t)|^{p}\right)dx
\end{equation}
and
\begin{equation}\label{Qdef}
\mathcal{Q}(u_{j})=\int_{-\infty }^{\infty }|u_{j}(x,t)|
^{2}\ dx, \ 1\leq j\leq 3.
\end{equation}

\smallskip

\noindent The mathematically exact theory for the nonlinear stability of solitary-wave solutions
began with Benjamin's theory \cite{[Benj]} (see also Bona \cite{[Bona]}) for the Korteweg-de
Vries equation.
After their papers on KdV and the regularized long-wave equations,
there are numerous literatures that have been devoted to the study of stability of solitary-wave
solutions for a variety of nonlinear dispersive equations. In particular,
Cazenave and Lions \cite{[CL]} developed an alternate approach to proving the stability of solitary waves
when they
are minimizers of the energy functional and when a compactness condition on minimizing sequences holds.
Their approach makes
use of the concentration compactness principle of Lions \cite{[L]} and
has the advantage of requiring less detailed analysis than the local
methods.
The Cazenave and Lions method has since been adapted by many different authors to prove the stability
results of a variety of nonlinear dispersive and wave equations
(see, for example, \cite{[AAu], [AB11], [San4], [San3], [Ngu], [Ngu2], [Ohta]} and references therein).

\smallskip

\noindent We now summarize the known results on the stability of solitary-wave
solutions of the coupled NLS systems. First, we provide some important results concerning two-component
NLS solitary waves that are relevant in our work.
In the case when $p=2, a_{11}=a_{22}=1,$ and $a_{21}=a_{12}=\beta >-1,$ the system \eqref{ODE} is known to have
explicit semi-trivial solution $(\Phi_1,\Phi_2,0)$ of the form
\begin{equation}\label{solformu}
\Phi_1(x)=\Phi_2(x)=\Phi_{\Omega}(x)=\sqrt{\frac{2\Omega}{1+\beta}} \textrm{sech}(\sqrt{\Omega}x),\ \Omega>0.
\end{equation}
In \cite{[Ohta]}, Ohta proved a stability result for two-component NLS solitary
waves of the form $e^{i\Omega t}\Phi_{\Omega}(x) \left(1,1 \right),$ for some $\Omega>0.$
Notice that since \eqref{CNLS} is invariant under the Galilean transformations
\begin{equation*}
u_{j}(x,t)\mapsto e^{-i\sigma ^{2}t+idx}u_{j}(x-2\sigma t,t),\ d,\sigma
\in \mathbb{R},\ j=1,2,3,
\end{equation*}
and the phase transformations
$u_{j}(x,t)\mapsto e^{i\beta _{j}}u_{j}(x,t),\ \beta _{j}\in \mathbb{R},$
one can write the solitary wave $e^{i\Omega t}\Phi_{\Omega}(x) \left(1,1 \right)$ into the following form
\begin{equation*}
\left(u_{1},u_{2}\right)=e^{i(\Omega-\sigma ^{2})t+i\sigma x}\Phi_{\Omega} (x-2\sigma t) \left(e^{i\beta_1},
e^{i\beta_2} \right).
\end{equation*}
When $a_{21}=a_{12}=\beta>0,\ \beta <\min\{a_{11},a_{22}\}$ or $\beta >\max\{a_{11},a_{22}\}$
and $\beta^{2}>a_{11}a_{22},$ Nguyen and Wang \cite{[Ngu]} proved the stability of solutions of the form
\begin{equation*}
e^{i\Omega t}\phi_{\Omega}(x) \left(\sqrt{\frac{\beta - a_{22}}{\beta^{2}-a_{11}a_{22}}},
\sqrt{\frac{\beta - a_{22}}{\beta^{2}-a_{11}a_{22}}}\right),
\end{equation*}
where $\phi_{\Omega}(x)$ is defined as in \eqref{solformu} with $\beta=0.$
We also note that the same authors (see \cite{[Ngu9]}) have
proved the stability of a two-parameter family of solitary waves for
two-component version of \eqref{CNLS}
in the special case $p=2,$ using the same method as in \cite{[AB11]}.
Similar techniques
have been used in \cite{[San3]} to prove the stability of (positive) ground-state solutions
of a more general two-component coupled NLS equations with power-type nonlinearities.

\smallskip

\noindent For 3-coupled systems such as \eqref{CNLS}, there are a variety of interesting results concerning the
existence of nontrivial solutions.
However, to our knowledge, the only available works regarding
the stability of nontrivial solutions for the full systems of type
\eqref{CNLS} are the papers \cite{[Ngu2],[Ngu5]}.
In \cite{[Ngu2]}, Nguyen and Wang considered \eqref{CNLS} in the special case when $p=2,$ and
proved the stability (in the sense defined above) of solutions, given by
\begin{equation}\label{NguSO2}
\sqrt{2\Omega}\ e^{i\Omega t}\textrm{sech}(\sqrt{\Omega}x)(\alpha_1, \alpha_2, \alpha_3), \ \Omega>0,
\end{equation}
under certain conditions
on coefficients $a_{jk}.$
More precisely, they made the following assumptions on
the matrix $B=(a_{jk})$ of positive coefficients:
\begin{enumerate}
\item[(1)] $B$ is invertible and the linear system $B\vec{\alpha}=\vec{1}$ is
solvable for $\vec{\alpha}=(\alpha_{1}^{2}, \alpha_{2}^{2},\alpha_{3}^{2}) \in \mathbb{R}_{+}^{3},$ where $\vec{1}=(1,1,1);$

\smallskip

\item[(2)] For all pair $j \neq k,$ $a_{jk}<\min\{a_{jj}, a_{kk}\}$ and $\det B_{jk}$ has the sign of $(-1)^{j+k+1}$ for $j \neq k ;$

\smallskip

\item[(3)] For all pair $j \neq k,$ $a_{jk}>\max\{a_{jj}, a_{kk}\}$ and
$\det B_{jk}$ has the sign of $(-1)^{j+k}$ for $j \neq k.$
\end{enumerate}
Then, using Lions' concentration compactness principle, they proved
that if the matrix $B$ satisfies (1) and (2) or (3),
the solutions of \eqref{CNLS} of the form \eqref{NguSO2}
are stable in $Y.$
The method used by Nguyen and Wang in \cite{[Ngu2]}
uses techniques from \cite{[Ngu]} with crucial ideas that
the constraints on the $L^2$ norms of components are not independently prescribed
and that the matrix of coefficients $B$ gives rise to positive numbers $\alpha_j$ such that the Euler-Lagrange equations can be
rewritten as uncoupled equations.

\smallskip

\noindent The paper \cite{[Ngu5]} is concerned with the stability of certain form of travelling-wave
solutions to $m$-component version of \eqref{CNLS} with $a_{ii}=a, a_{kj}=b,$ and $a+2b>0.$
Their results generalize the ones obtained in \cite{[Ngu], [Ngu2]} to include a more general case of coupled
nonlinear Schr\"{o}dinger equations. To state the precise statement of their stability result,
for any $\omega_1=\omega_2=\Omega>0,$ set
\begin{equation*}
\phi_{\Omega, a+2b}(x)=\left(\frac{\Omega}{a+2b}\right)^{1/(2p-2)}\phi(\sqrt{\Omega}\ x),
\end{equation*}
where $\phi(x)$ is the unique positive, spherically symmetric, and
decreasing solution of
\begin{equation*}
\left\{
\begin{aligned}
- & u_{xx}+u-|u|^{2p-2}u=0, \ x\in \mathbb{R}, \\
& u \in H^1(\mathbb{R})\setminus \{0\}.
\end{aligned}
\right.
\end{equation*}
It has been shown in
\cite{[Ngu5]} that when $b>0$ and $a+2b>0,$ and for $2 \leq p < 3,$
travelling-wave solutions
to \eqref{CNLS} of the form $e^{i\Omega t}\phi_{\Omega, a+2b}(x)(1, 1, 1)$
are stable in the following sense: for any $\epsilon>0,$ there
exists $\delta>0$ such that if $\mathbf{u}_{1,0}\in Y$ with
\begin{equation*}
\inf_{\gamma_j\in \mathbb{R}} \inf_{y \in \mathbb{R}} \left\{ \sum_{j=1}^{3}\|u_{j,0}
-e^{i\gamma_j}\phi_{\Omega, a+2b}( \cdot + y)\|_{1}\right\} < \delta,
\end{equation*}
then the solution $\mathbf{u}(x,t)$ of \eqref{CNLS} with initial
condition $\mathbf{u}(\cdot, 0)=\mathbf{u}_{1,0}$ satisfies
\begin{equation*}
\inf_{\gamma_j\in \mathbb{R}} \inf_{y \in \mathbb{R}} \left\{ \sum_{j=1}^{3}\|u_{j}
-e^{i\theta_j}\phi_{\Omega, a+2b}( \cdot + y)\|_{1}\right\} < \epsilon,
\end{equation*}
uniformly for all $t\geq 0.$

\smallskip

\noindent As seen from the preceding discussion, the stability results obtained from all these papers \cite{[Ngu], [Ngu2], [Ngu5],[Ohta]} are for
one-parameter family of solitary waves, in which each component is a multiple
of a hyperbolic secant function. Their stability results were obtained by characterizing solitary
waves as minimizers of an energy functional subject to the constraints that were not independently chosen.
In this paper we study the variational problem of finding the extremum of
the energy functional $\mathcal{H}(u_1,u_2,u_3)$ satisfying the constraints
\begin{equation*}
\int_{-\infty}^{\infty}|u_1|^2\ dx =r,\ \int_{-\infty}^{\infty}|u_2|^2\ dx =s,\ \textrm{and}\ \int_{-\infty}^{\infty}|u_3|^2\ dx =t.
\end{equation*}
Such solutions are of interest in physics and sometimes referred to as normalized solutions and the associated
solitary waves as normalized solitary waves.
Our method leads to the existence and stability results concerning a fully nontrivial
three-parameter family of solitary waves.
To the best
of our knowledge, the results of this paper are the first existence and stability results for
such normalized solitary waves of three-component nonlinear systems.
The reader may see \cite{[AAu],[AB11], [San3], [Ngu9]} for the existence and stability results of
independently prescribed $L^2$-norm solutions
to two-component systems.

\smallskip

\noindent We now describe the main results of this paper. We prove that
the full equations \eqref{CNLS} has
a non-empty stable set of positive normalized solitary-wave solutions for all positive constants $a_{kj}=a_{jk}$ and all $p\in[2,3)$
(we say that a solution of \eqref{CNLS} is positive
if each component is in the form $e^{i\theta}\varphi(x),$
where $\theta$ is a real constant and $\varphi$ is an $\mathbb{R}$-valued
positive function).
The existence result is obtained via a variational approach
and using Cazenave-Lions method \cite{[L], [CL]}.
The parameters $\omega_1,\omega_2,\omega_3\in \mathbb{R}$ in the equations \eqref{ODE} appear as
Lagrange multipliers.
More precisely, let $\mathcal{H}$ and $\mathcal{Q}$ be as defined in \eqref{Edef} and \eqref{Qdef}, respectively;
it is easy to see using the Sobolev embedding theorem that $\mathcal{H}$ and $\mathcal{Q}$ define
continuous maps from $Y$ to $\mathbb{R}.$
For $r,s,t>0,$ let
\begin{equation}
\Delta_{r,s,t}=\left\{\mathbf{f} \in Y: \|f_1\|^2=r,\ \|f_2\|^2=s,\ \text{and} \ \|f_3\|^2=t \right\}
\end{equation}
and define the function $\lambda(r,s,t)$ by
\begin{equation}
\label{Ist}
\lambda(r,s,t)=\inf \left\{ \mathcal{H}(\mathbf{f}): \mathbf{f} \in \Delta_{r,s,t} \right\}.
\end{equation}
A minimizing sequence for $\lambda(r,s,t)$ is any sequence
$\{\mathbf{f}_n\}$ in $Y$ satisfying the conditions
\begin{equation*}
\lim_{n\to \infty}\|f_{1, n}\|^2=r,\ \lim_{n\to \infty}\|f_{2, n}\|^2=s,\ \lim_{n\to \infty}\|f_{3, n}\|^2=t,\ \textrm{and}\ \lim_{n\to \infty}\mathcal{H}(\mathbf{f}_n)=\lambda(r,s,t).
\end{equation*}
To each minimizing sequence $\{\mathbf{f}_{n}\}$ of the problem $\lambda(r,s,t),$  we associate a
sequence of nondecreasing functions
$P_{n}:[0,\infty)\to [0,r+s+t]$ defined by
\begin{equation*}
P_n(\eta)=\sup_{y\in \mathbb{R}}\int_{y-\eta}^{y+\eta}\rho _{n}(x)\ dx,
\end{equation*}
where $\rho _{n}(x)$ is given by
\begin{equation*}
 \rho _{n}(x):=|f_{1,n}(x)|^2 + |f_{2,n}(x)|^2+|f_{3,n}(x)|^2.
 \end{equation*}
A standard argument shows that any uniformly bounded sequence of nondecreasing functions on
$[0,\infty)$ must have a subsequence which converges pointwise
to a nondecreasing limit function on $[0,\infty).$
Hence $P_{n}(\eta)$
has such a subsequence, which we again denote by $P_n.$
Let $P(\eta):[0,\infty)\to [0,r+s+t]$ be the nondecreasing function to which $P_n$ converges, and define
\begin{equation}
\label{defgamma} \gamma =\lim_{\eta\to \infty }P(\eta).
\end{equation}
Then $\gamma$ satisfies $0\leq \gamma \leq r+s+t.$
The method of Cazenave and Lions \cite{[CL]}, as applied to this situation,
consists of two observations. The first is that if $\gamma = r+s+t,$
then the minimizing sequence $\{\mathbf{f}_{n}\}$ has a
subsequence which, when its terms are suitably translated, converges
strongly in $Y$ to an element of the set $\mathcal{O}_{r,s,t}$ defined by
\begin{equation}
\mathcal{O}_{r,s,t}=\left\{\mathbf{\Phi} \in Y: \mathcal{H}(\mathbf{\Phi})=\lambda(r,s,t), \mathbf{\Phi}\in \Delta_{r,s,t} \right\}.
\end{equation}
The second is that certain properties of
the variational problem imply that $\gamma$ must equal $r+s+t$ for every minimizing sequence
$\{\mathbf{f}_{n}\}.$ It follows that not only do minimizers exist in $Y,$ but every minimizing
sequence converges in $Y$ norm to the set $\mathcal{O}_{r,s,t}.$
Typically, one proves $\gamma=r+s+t$ by ruling out other two possibilities, namely $\gamma=0$ and $0<\gamma<r+s+t.$
A lemma (Lemma~2.10 of \cite{[AB11]}) concerning the symmetric rearrangement of functions plays an important role
in our proof. In Sections~\ref{SecVar} and \ref{SecVar1}, we provide the details
of the method.

\smallskip

\noindent We prove below (see Theorem~\ref{Main1}) that the problem \eqref{Ist} has a
solution in $\Delta_{r,s,t}$ for the range $2\leq p <3.$
In particular,
the set $\mathcal{O}_{r,s,t}$ is nonempty.
The set $\mathcal{O}_{r,s,t}$ consists of solitary-wave profiles for \eqref{CNLS}. More precisely,
if $(\Phi_1,\Phi_2,\Phi_3) \in \mathcal{O}_{r,s,t},$ they satisfies the Euler-Lagrange equations
\begin{equation*}
-\Phi_{j, xx}+\omega_j \Phi_j=\left(\sum_{k=1}^{3}a_{kj}|\Phi_k|^p\right)|\Phi_j|^{p-2} \Phi_j, \ j=1,2,3,
\end{equation*}
where $\omega_j$ are the Lagrange multipliers. The preceding equations are satisfied by $(\Phi_1,\Phi_2,\Phi_3)$
if and only if the triple
\begin{equation}\label{SOFinal}
(u_1(x,t), u_2(x,t), u_3(x,t))=(e^{i\omega_1 t}\Phi_1(x),e^{i\omega_{2} t}\Phi_{2}(x),e^{i\omega_{3} t}\Phi_{3}(x))
\end{equation}
is a solutions of \eqref{CNLS}, and since \eqref{CNLS} is invariant under the
Galilean transformations and the phase transformations, one can always write \eqref{SOFinal} into the form \eqref{SO}.

\smallskip

\noindent The question
about the characterization of the set $\mathcal{O}_{r,s,t}$ is addressed in
Section~\ref{SecVar1} (see Theorem~\ref{Main2}).
Namely, we prove that for each $\mathbf{\Phi} \in \mathcal{O}_{r,s,t}$ there
exists positive real-valued functions
$\phi_1, \phi_2, \phi_3 \in H^{1}$ such that
\begin{equation*}
\Phi_j(x)=\lambda_j \phi_j(x), \ \lambda_j \in S^1, \ j=1,2,3.
\end{equation*}
Also, the functions $\Phi_j$ are infinitely differentiable on $\mathbb{R}.$

\noindent Finally, Theorem~\ref{Main3} proves that
the $\mathcal{O}_{r,s,t}$ forms a stable set for the associated initial-value problem to \eqref{CNLS}
in the sense of Definition~\ref{def1}.

\medskip

\noindent {\bf Notation.} For $1\leq p\leq \infty ,$ the space of all complex-valued measurable functions whose $p-$th
power is integrable will be denoted by
$L^{p}=L^{p}(\mathbb{R})$ and its norm $\left\vert f\right\vert _{p}$ is given by
\begin{equation*}
\left\vert f\right\vert _{p}=\left( \int_{-\infty }^{\infty }\left\vert
f\right\vert ^{p}dx\right) ^{1/p}\textrm{ \ for }1\leq p<\infty,
\end{equation*}
and $\left\vert f\right\vert _{\infty }$ is the essential supremum of $%
\left\vert f\right\vert $ on $\mathbb{R}.$
Whether we intend the functions in $L^{p}$ to be real-valued or complex-valued will be clear from the context.
We denote
by $H^1=H^{1}(\mathbb{R})$ the Sobolev space of all complex-valued, measurable functions
on $\mathbb{R}$ such that both $f$ and
$f^{\prime}$ are in $L^{2}.$ The norm $\|.\|_{1}$ on the space $H^{1}$ is given by
 \begin{equation*}
 \|f\|_{1}=\left(\int_{-\infty }^{\infty }(|f(x)|^{2}+|f^{\prime}(x)|^{2}) \ dx \right)^{1/2}.
 \end{equation*}
In particular, we use $\left\Vert
f\right\Vert $ to denote the $L^{2}$ norm of a
function $f.$
We define $X$ to be the product space $X=H^{1}(\mathbb{R})\times H^{1}(
\mathbb{R})$ furnished with the product norm
\begin{equation*}
\|(f,g)\|_X^{2}=\int_{-\infty}^{\infty}\left(|f(x)|^{2}+|g(x)|^{2}\right) \ dx + \int_{-\infty}^{\infty}\left(|f^{\prime}(x)|^{2}+|g^{\prime}(x)|^{2}\right) \ dx
\end{equation*}
and the space $Y$ to be the space
$Y=H^{1}(\mathbb{R})\times H^{1}(
\mathbb{R}) \times H^{1}(\mathbb{R})$
equipped with the product norm, which we denote by $\|\cdot\|_Y.$ That is,
\begin{equation*}
\|(f_1,f_2,f_3) \|_Y^{2}
 = \sum_{j=1}^{3}\int_{-\infty}^{\infty}|f_j(x)|^{2} \ dx +\sum_{j=1}^{3}\int_{-\infty}^{\infty}|f^{\prime}_j(x)|^{2} \ dx.
\end{equation*}
If $T>0$ and $Z$ is any Banach space, we denote by $\mathcal{C}([0,T]; \ Z)$ the Banach space
of continuous maps $f:[0,T] \to Z,$ with norms given by
\begin{equation*}
\|f\|_{\mathcal{C}([0,T]; \ Z)} = \sup_{t\in [0,T]}\|f(t)\|_Y.
\end{equation*}
For notational convenience, we set $\mathbf{f}=(f_1,f_2,f_3),\ \mathbf{f}_{n}=(f_{1,n},f_{2,n}, f_{3,n}),$
and $\mathbf{\Phi} =(\Phi_1,\Phi_2,\Phi_3).$
We denote by $S^1$ the set of all complex numbers of the form $e^{i\theta},$ i.e.,
\begin{equation*}
S^1=\{e^{i\theta} : \theta \in \mathbb{R}\}.
\end{equation*}
 The letter $C$ will be used to denote various
constants whose exact values are immaterial and which may vary from one line to the
next.

%The title of your section 2
\section{The Variational Problem}\label{SecVar}

\noindent Throughout this section, we assume that $p \in [2,3)$ and $a_{jk}$ are positive real constants
satisfying $a_{kj}=a_{jk}$ for all $k,j\in\{1,2,3\}.$

\smallskip

\noindent We first establish some properties of the variational problem and its minimizing
sequences which are independent of the value $\gamma.$

\begin{lem} \label{Ibounded} If $\{\mathbf{f}_n\}$ is a minimizing
sequence for $\lambda(r,s,t),$ then
there exists a constant $B$ such that
\begin{equation*}
\sum_{j=1}^{3} \|f_{j,n}\|_1 \leq B.
\end{equation*}
Moreover, for any $r, s, t>0,$ one has
$-\infty <\lambda(r,s,t)<0.$
\end{lem}
\noindent{\bf Proof.}
Let $ \{\mathbf{f}_{n}\}$ be a minimizing
sequence for the problem $\lambda(r,s,t)$. Then $\|f_{1,n}\|, \|f_{2,n}\|,$ and $\|f_{3,n}\|$ are
bounded. Using the Gagliardo-Nirenberg inequality,  we obtain
\begin{equation}
\label{fLqbound} \int_{-\infty}^{\infty}|f_{j,n}|^{2p} \ dx \leq
C\left( \int_{-\infty}^{\infty}|f_{j,nx}|^{2} \ dx\right)^{(p-1)/2}\left( \int_{-\infty}^{\infty} |f_{j,n}|^2 \ dx \right)^{(p+1)/2},
\end{equation}
where $C=C(p,r,s,t).$ For each $j,k= 1, 2, 3,$ using the Cauchy-Schwarz inequality, we also have
\begin{equation}\label{mixtermbound}
 \int_{-\infty }^{\infty }|f_{j, n}|^{p}|f_{k,n}|^{p}\ dx
\leq \frac{1}{2}\left(\int_{-\infty }^{\infty } |f_{j, n}|^{2p} \ dx + \int_{-\infty }^{\infty }|f_{k,n}|^{2p} \ dx \right).
\end{equation}
Now
\begin{equation*}
\begin{aligned}
& \|\mathbf{f}_{n}\|_Y^2 =\|f_{1,n}\|_1^2+\|f_{2,n}\|_1^2+\|f_{3,n}\|_1^2 \\
&\ \ =\mathcal{H}(\mathbf{f}_{n})+\frac{1}{p}\int_{-\infty }^{\infty
}\sum_{k,j=1}^{3}a_{kj}|f_{k,n}|^{p}|
f_{j,n}|^{p}+(r+s+t).
\end{aligned}
\end{equation*}
Since $ \{\mathbf{f}_{n}\}$ is a
minimizing sequence, $\mathcal{H}(\mathbf{f}_{n})$ is bounded.  Using \eqref{fLqbound} and \eqref{mixtermbound}, it follows that
\begin{equation*}
\|\mathbf{f}_{n}\|_Y^2 \le C\left(1 +
\|\mathbf{f}_{n}\|_Y^{p-1}\right),
\end{equation*}
and hence, the existence of the desired bound $B$ follows.

\noindent The claim that $\lambda(r,s,t)>-\infty$ easily follows using the estimates \eqref{fLqbound} and \eqref{mixtermbound}.
To prove $\lambda(r,s,t)<0,$ choose any function $f_1 \in H^{1}$ such that $%
\|f_1\|^{2}=r.$ Let $f_2=(s/r)^{1/2}f_1$ and $%
f_3=(t/r)^{1/2}f_1.$ Then
$\|f_2\|^{2}=s$ and $\|f_3\|
^{2}=t.$ For each $\theta >0,$ define
$f_{j, \theta }=\theta ^{1/2}
f_{j}(\theta x),\ j=1,2,3.$
Then, for all $%
\theta ,$ we have
\begin{equation*}
\|f_{1, \theta }\|^{2}=\|f_1\|
^{2}=r,\ \|f_{2, \theta }\|^{2}=\|
f_2\|^{2}=s,\ \text{and}\ \|f_{3, \theta }\|
^{2}=\|f_3\|^{2}=t,
\end{equation*}
and
\begin{equation*}
\mathcal{H}(\mathbf{f}_{\theta })=\int_{-\infty }^{\infty }\left( \theta
^{2}\sum_{j=1}^{3}|f_{j, x}|^{2}-\theta^{p-1}
\sum_{k,j=1}^{3}\frac{1}{p} \ a_{kj}|f_k|^{p}|
f_j|^{p}\right) dx.
\end{equation*}
Substituting $f_2$ and $f_3$ as defined above on the right hand side, we obtain
\begin{equation*}
\begin{aligned}
\mathcal{H}(\mathbf{f}_{\theta })&=\int_{-\infty }^{\infty }\left[ \theta ^{2}\left( 1+\frac{s%
}{r}+\frac{t}{r}\right) |f_{1, x}|^{2}\right.
-\frac{1}{p} \ \theta^{p-1} \left( a_{11}+a_{22}\frac{s^p}{r^p} \right. \\
\ & +a_{33}\frac{t^p}{r^p}
\ +\left. \left. 2a_{12}\frac{s^{p/2}}{r^{p/2}}+2a_{13}\frac{t^{p/2}}{r^{p/2}}+ 2a_{23}\frac{s^{p/2}}{r^{p/2}}
\frac{t^{p/2}}{r^{p/2}} \right) |f_1|^{2p}\right] \ dx
\end{aligned}
\end{equation*}
From our assumption on the constants $a_{kj},$ the quantity
\begin{equation*}
a_{11}+a_{22}\frac{s^p}{r^p}
+a_{33}\frac{t^p}{r^p}
\ +2a_{12}\frac{s^{p/2}}{r^{p/2}}+2a_{13}\frac{t^{p/2}}{r^{p/2}}+ 2a_{23}\frac{s^{p/2}}{r^{p/2}}
\frac{t^{p/2}}{r^{p/2}} > 0.
\end{equation*}
Hence, we can get $\mathcal{H}(\mathbf{f}_{\theta_0 })<0$
by choosing $\theta =\theta_{0}$ sufficiently small. \hfill $\Box$

\medskip

\begin{lem} \label{bdgnxbelow}
For any minimizing sequence $ \{\mathbf{f}_{n}\}\subset Y$
for $\lambda(r,s,t)$, the following statements hold for all large enough $n,$

\noindent{(i) if $r>0$ and $s, t \ge 0$,
then there exists $\delta_{1}>0$ such that $\|f_{1, nx}\|\geq \delta_{1}.$}

\noindent{(ii) if $s>0$ and $r, t \ge 0$,
then there exists $\delta_{2}>0$ such that $\|f_{2, nx}\|\geq \delta_{2}.$}

\noindent{(iii) if $t>0$ and $r, s \ge 0$,
then there exists $\delta_{3}>0$ such that $\|f_{3, nx}\|\geq \delta_{3}.$}
\end{lem}
\noindent{\bf Proof.}
To prove statement (i), suppose to the contrary that no such
constant $\delta_1$ exists. Then, by taking a subsequence if necessary, one may assume that
$\displaystyle \lim_{n \to \infty} \|f_{1, nx}\| = 0$.
By the Gagliardo-Nirenberg inequalities,
\begin{equation*}
\lim_{n \to \infty} \int_{-\infty }^{\infty }|f_{1, n}|^{p}|f_{j, n}|^{p}\ dx=0, \ \ j=1,2,3.
\end{equation*}
Therefore, we have that
\begin{equation}
\label{Ibbelow}
\lambda(r,s,t) =\lim_{n\to \infty }\mathcal{H}(\mathbf{f}_n)
=\lim_{n\to \infty }\left(\sum_{j=2}^{3} \|f_{j, nx}\|^{2}-\frac{1}{p}\sum_{k, j =2}^{3}a_{kj}|f_{j, n}f_{k, n}|_{p}^{p}\right).
\end{equation}
Pick any $\psi_1\geq 0$ such that $\|\psi_1\|^2=r.$ For every $\theta >0$,
define $\psi_{1, \theta}(x)=\theta^{1/2}\psi_{1}(\theta x).$
Then $\|\psi_{1, \theta} \|^2=r,$ and hence, for all $n,$
\begin{equation*}
\lambda(r,s,t) \le \mathcal{H}(\psi_{1, \theta}, f_{2, n}, f_{3, n}).
\end{equation*}
On the other hand, if one defines
\begin{equation}\label{elzero}
\eta = \theta^2\int_{-\infty }^{\infty }|\psi_{1, x}|^2\ dx-
\theta^{p-1} \int_{-\infty }^{\infty }\frac{a_{11}}{p}|\psi_1|^{2p} \ dx,
\end{equation}
then $\eta < 0$ for sufficiently small $\theta.$
With this notation, we have for all $n\in\mathbb{N}$,
\begin{equation*}
\begin{aligned}
\lambda(r,s,t) &\leq \mathcal{H}(\psi_{1, \theta}, f_{2, n}, f_{3, n}) \\
& \le \int_{-\infty}^{\infty } \left(\sum_{j=2}^{3}|f_{j, nx}|^{2}-\frac{1}{p}\sum_{k, j =2}^{3}a_{kj}|f_{j, n}|^{p}
|f_{k, n}|^{p}
\right) \ dx + \eta.
\end{aligned}
\end{equation*}
Consequently
\begin{equation*}
\begin{aligned}
\lambda(r,s,t) &\leq \lim_{n \to \infty}\int_{-\infty}^{\infty
}\left(\sum_{j=2}^{3}|f_{j, nx}|^{2}-\frac{1}{p}\sum_{k, j =2}^{3}a_{kj}|f_{j, n}|^{p}
|f_{k, n}|^{p}
\right) \ dx
+ \eta,
\end{aligned}
\end{equation*}
which contradicts \eqref{Ibbelow} and \eqref{elzero}, and hence, statement (i) follows.
The statements (ii) and (iii) can be proved similarly. \hfill $\Box$

\medskip

\begin{lem}\label{minforJ}
\label{Is0} Let $\alpha_1, \alpha_2 , \beta > 0$ and $p\in [2,3).$ For $f,g\in H^1(\mathbb{R}),$ define the functional
\begin{equation}\label{defJ}
F(f,g)=\|f_x\|^2 +\|g_x\|^2 -\frac{1}{p}\left(\alpha_1
|f|_{2p}^{2p}+ \alpha_2|g|_{2p}^{2p}+ 2\beta|fg|_{p}^{p} \right).
\end{equation}
Then for any $a_1, a_2 >0$ there exists a nontrivial solution to the problem
\begin{equation}\label{mdefm}
m(a_1,a_2)=\inf\{F(f,g):\|f\|^2=a_1\ \textrm{and}\ \|g\|^2=a_2\}.
\end{equation}
Furthermore, if $(\tilde{\phi}_{a_1},\tilde{\phi}_{a_2})$ is a
solution of \eqref{mdefm}, then there exists $\theta_1,\theta_2\in \mathbb{R}$ and
positive real-valued functions $\phi_{a_1}$ and $\phi_{a_2}$ such that
$\tilde{\phi}_{a_1}(x)=e^{i\theta_1}\phi_{a_1}(x)$ and $\tilde{\phi}_{a_2}(x)=e^{i\theta_2}\phi_{a_2}(x).$
In particular,
\begin{equation}\label{twopos}
F(\phi_{a_1}, \phi_{a_2})=m(a_1,a_2).
\end{equation}
\end{lem}
\noindent{\bf Proof.} Let $\{(f_n, g_n) \}$ be a sequence of functions in $X$ satisfying the conditions
$\underset{n\to \infty}\lim \|f_n\|^2 =a_1, \underset{n\to \infty}\lim\|g_n\|^2= a_2,$
 and
\begin{equation*}
\lim_{n \to \infty} F(f_n, g_n) = m(a_1,a_2).
\end{equation*}
Using the concentration compactness argument,
it has been proved in \cite{[San3]} that the sequence
$(f_n,g_n)$ converges,
up to taking a subsequence and after suitable translations, strongly to some function $(\tilde{\phi}_{a_1},\tilde{\phi}_{a_2})$ in $X$ norm.
Then the pair $(\tilde{\phi}_{a_1},\tilde{\phi}_{a_2})$ achieves the minimum,
\begin{equation}
\label{E220}
F(\tilde{\phi}_{a_1},\tilde{\phi}_{a_2})=m(a_1,a_2)
\end{equation}
and must satisfy the Euler-Lagrange equation
\begin{equation*}
\left\{
\begin{aligned}
-\tilde{\phi}^{\prime \prime}_{a_1}+\mu_1\tilde{\phi}_{a_1} & = \alpha_1 |\tilde{\phi}_{a_1}|^{2p-2}\tilde{\phi}_{a_1}+\beta |\tilde{\phi}_{a_2}|^{p}\tilde{\phi}_{a_1}, \\
-\tilde{\phi}^{\prime \prime}_{a_2}+\mu_2\tilde{\phi}_{a_2} & = \alpha_2|\tilde{\phi}_{a_2}|^{2p-2}\tilde{\phi}_{a_2}+\beta |\tilde{\phi}_{a_1}|^{p}\tilde{\phi}_{a_2},
\end{aligned}
\right.
\end{equation*}
for some real numbers $\mu_2$ and $\mu_2.$ As in the proof of Theorem~\ref{Main2} below, one can show that
there exists real numbers $\theta_1, \theta_2,$ and real-valued positive
functions $\phi_{a_1}$ and $\phi_{a_2}$ such that
$\tilde{\phi}_{a_1}(x)= e^{i\theta_1}\phi_{a_1}(x)$ and $\tilde{\phi}_{a_2}(x)=e^{i\theta_2}\phi_{a_2}(x).$
Then, for some sequence $\{y_n\}$ of real numbers,
\begin{equation*}
\{(e^{-i\theta_1}f_{n_k}(\cdot + y_k),e^{-i\theta_2}g_{n_k}(\cdot + y_k)) \}
\end{equation*}
converges
strongly to $(\phi_{a_1}, \phi_{a_2})$ in $X.$
Finally, since
$F(\phi_{a_1}, \phi_{a_2})=F(\tilde{\phi}_{a_1},\tilde{\phi}_{a_2}),$ the assertion \eqref{twopos} follows
from \eqref{E220}. \hfill $\Box$

\medskip

\begin{lem} \label{negative}
For $r, s, t>0,$ let $\{\mathbf{f}_n\}\subset Y$ be a minimizing sequence for the problem $\lambda(r,s,t)$.
Then, for each $j=1,2,3,$ there exists $\delta_{j}>0$ such that for all large enough
$n$,
\begin{equation*}
\int_{-\infty}^\infty\left(|f_{j, nx}|^2 -\frac{1}{p}|f_{j, n}|^{p}\sum_{k =1}^{3}a_{jk}|f_{k, n}|^{p}
\right)\ dx \le -\delta_{j}.
\end{equation*}
\end{lem}
\noindent{\bf Proof.} We prove the lemma for $j=1.$ The proofs for cases $j=2,3$ are similar.
Suppose the conclusion of lemma is false. Then, by taking a 
subsequence if necessary, one can find a minimizing sequence
$\{\mathbf{f}_n\}$ for the problem $\lambda(r,s,t)$ such that
\begin{equation}
\liminf_{n \to \infty} \int_{-\infty}^\infty\left(|f_{1, nx}|^2 -\frac{1}{p}|f_{1, n}|^{p}\sum_{k =1}^{3}a_{1k}|f_{k, n}|^{p}
\right)\ dx \ge 0,
\end{equation}
and hence
\begin{equation}
\begin{aligned}
\lambda(r,s,t) &= \lim_{n\to\infty}\mathcal{H}(f_{1,n},f_{2,n}, f_{3,n}) \\
&\ge \liminf_{n \to \infty}
\int_{-\infty}^\infty \left(\sum_{j=2}^{3}|f_{j, nx}|^{2}-\frac{1}{p}\sum_{k, j =2}^{3}a_{kj}|f_{j, n}|^{p}
|f_{k, n}|^{p}
\right) \ dx
\end{aligned}
 \label{IgeJg0}
\end{equation}
\noindent Let $F, \phi_s,$ and $\phi_t$ be as defined in Lemma \ref{minforJ} with $f=f_{2,n}, g=f_{3,n},
\alpha_1=a_{22}, \alpha_2 =a_{33}, \beta =a_{23}, a_1=s,$ and $a_2=t.$  Then
\eqref{IgeJg0} gives that
\begin{equation}
\lambda(r,s,t) \ge F(\phi_s, \phi_t). \label{IgeJg02}
\end{equation} On the other hand, take any $f_1 \in H^1$ such that $\|f_1 \|^2=r$ and
\begin{equation}\label{minf2}
\int_{-\infty}^{\infty}\left(|f_{1, x}|^2 -\frac{a_{11}}{p}|f_1|^{2p}-\frac{2a_{12}}{p} |f_1|^p | \phi_s|^{p}
-\frac{2a_{13}}{p}|f_1|^{p} |\phi_t|^{p} \right)\ dx < 0.
\end{equation}
To construct such a function $f_1,$ take an arbitrary smooth function $\psi\geq 0$ with compact support
which satisfies $\psi(0)=1$ and $\|\psi\|=r,$ and for $\theta > 0,$ define $\psi_{\theta}(x)=\theta^{1/2}\psi(\theta x).$
Then, the function $f_1=\psi_{\theta}$ satisfies \eqref{minf2} for sufficiently small $\theta.$
With the use of \eqref{minf2}, we have that
\begin{equation*}
\begin{aligned}
& \lambda(r,s,t) \le \mathcal{H}(f_1,\phi_s, \phi_t) \\
 & \le \int_{-\infty}^{\infty}\left(|f_{1, x}|^2 -\frac{a_{11}}{p}|f_1|^{2p}-\frac{2a_{12}}{p} |f_1|^p | \phi_s|^{p}
-\frac{2a_{13}}{p}|f_1|^{p} |\phi_t|^{p} \right)\ dx + F(\phi_s, \phi_t) \\
 & < F(\phi_s, \phi_t),
\end{aligned}
\end{equation*}
which contradicts \eqref{IgeJg02} and hence, lemma follows. \hfill $\Box$

\medskip

\noindent In what follows we denote by $w^\ast$ the symmetric rearrangement
of a measurable function $w:\mathbb{R} \to [0,\infty)$ (for the definition and a
concise exposition of the basic properties of symmetric
rearrangements, we refer the reader to the excellent book by Lieb and Loss \cite{[LLieb]}).
 Notice that if
$\mathbf{f}$ belongs to $Y$, then all $|f_j|$ also belong to $H^1$, and hence
the rearrangements $|f_j|^\ast$ are well-defined.
We also note the following basic property about rearrangements
\begin{equation}\label{rearrangement1}
\int_{-\infty }^{\infty }(|f_j|^\ast)^{2p}\ dx = \int_{-\infty }^{\infty }|f_j|^{2p}\ dx.
\end{equation}

\begin{lem} \label{posmin}
The following hold for all $(f_1,f_2,f_3) \in Y,$
\begin{equation}\label{Hdecreasing}
\mathcal{H}(|f_1|,|f_2|,|f_3|)\leq \mathcal{H}(f_1,f_2,f_3)
\end{equation}
and
\begin{equation}\label{Hdecreasing2}
\mathcal{H}(|f_1|^\ast,|f_2|^\ast, |f_3|^\ast)\leq \mathcal{H}(f_1,f_2,f_3).
\end{equation}
\end{lem}
\noindent{\bf Proof.}
The proof of \eqref{Hdecreasing} follows from a basic fact of analysis that
\begin{equation}\label{decrea2}
\int_{-\infty }^{\infty }\left\vert \left\vert g\right\vert _{x}\right\vert
^{2}\ dx\leq \int_{-\infty }^{\infty }\left\vert g_{x}\right\vert ^{2}\ dx.
\end{equation}
(For a proof of \eqref{decrea2}, see Lemma~3.5 of \cite{[ABS]}.) To prove \eqref{Hdecreasing2},
we use the following inequality about rearrangements
\begin{equation}\label{rearrangement2}
\int_{-\infty }^{\infty }(|f_j|^\ast) ^{p}(|f_k|^\ast)^{p} \ dx \geq
\int_{-\infty }^{\infty }|f_j|^{p}|f_k|^{p}\ dx,\ \forall k, j=1,2,3\ \textrm{and}\ k \neq j.
\end{equation}
We also
have that (see, for example, Lieb and Loss\cite{[LLieb]})
\begin{equation*}
\int_{-\infty }^{\infty }| (|f_j|^\ast)_{x}| ^{2}\ dx\leq
\int_{-\infty }^{\infty }||f_j|_{x}|^{2}\ dx.
\end{equation*}
Using all the facts above and \eqref{rearrangement1},
the inequality \eqref{Hdecreasing2} easily follows from \eqref{Hdecreasing}. \hfill $\Box$

\medskip

\noindent The next lemma is one-dimensional version of Proposition~1.4 of \cite{[Byeo]}.
A proof of this lemma is given in \cite{[AB11]}.

\begin{lem} Let $f$ and $g$ be the functions such that

\smallskip

(i) $f\geq 0$ and $g\geq 0$ for all $x\in \mathbb{R}.$

\smallskip
(ii) are even, $C_c^\infty,$ and non-increasing on $[0,\infty).$

\smallskip

\noindent Suppose $x_1$ and $x_2$ are the numbers such that the translated functions $f(x+x_1)$ and
$g(x+x_2)$ have disjoint supports, and define
$ w(x)=f(x+x_1)+g(x+x_2).$
Then $(w^\ast)'\in L^2$
and satisfies
\begin{equation}
\|(w^\ast)'\|^2 \le \|w'\|^2 - \frac34 \min\{\|f'\|^2,\|g'\|^2\},
\label{garineq}
\end{equation}
where the derivative is understood in the distribution sense.
\label{garlem}
\end{lem}

\noindent We now prove that the function $\lambda(r,s,t)$ is strictly subadditive:
\begin{lem} \label{subadd}
Let $r_1, r_2, s_1,s_2,t_1,t_2\geq 0$ be given, and suppose further that
$r_1+r_2>0, s_1+s_2>0$, $t_1+t_2>0$, $r_1+s_1+t_1>0$, and $r_2+s_2+t_2>0$. Then
\begin{equation}
\label{SUBA}
\lambda(r_{1}+r_{2},s_{1}+s_{2},t_{1}+t_{2})<\lambda(r_{1},s_{1},t_{1})+\lambda(r_2,s_{2},t_{2}).
\end{equation}
\end{lem}
\noindent{\bf Proof.} We follow closely the arguments used in \cite{[AB11]}. For $i=1,2,$ we first
construct minimizing sequences $(f_{1, n}^{(i)},f_{2, n}^{(i)},f_{3, n}^{(i)})$ for the problem $\lambda(r_i,s_i,t_i)$
such that for each $j \in \{1,2,3\}$ and all $n \in \mathbb{N},$ the functions $f_{j, n}^{(i)}$ are $\mathbb{R}$-valued,
satisfy conditions (i) and (ii) of Lemma~\ref{garlem},
and
\begin{equation*}
\|f_{1, n}^{(i)}\|^2=r_i,\ \|f_{2, n}^{(i)}\|^2=s_i, \ \text{and} \ \|f_{3, n}^{(i)}\|^2=t_i.
\end{equation*}
Without loss of generality, take $i=1$, since the case $i=2$ is exactly similar.
We may also assume that $r_1>0, s_1> 0,$ and $t_1> 0$,
otherwise we can just take $f_{1, n}^{(1)}, f_{2, n}^{(1)},$ or $f_{3, n}^{(1)}$ to be identically zero on $\mathbb{R}$ for all $n.$
Let $(e_n^{(1)}, p_n^{(1)},q_n^{(1)})$ be any minimizing sequence for $\lambda(r_1,s_1,t_1)$.
 By the continuity of $\mathcal{H}$ and the density of compactly supported functions in $H^1$,
 we can approximate $(e_n^{(1)}, p_n^{(1)},q_n^{(1)})$ by compactly supported functions
 $(e_n^{(2)}, p_n^{(2)},q_n^{(2)}).$ Then $(e_n^{(2)}, p_n^{(2)},q_n^{(2)})$ forms a
 minimizing sequence for the problem $\lambda(r_1,s_1,t_1).$
 Define
\begin{equation*}
(e_n^{(3)},p_n^{(3)},q_n^{(3)})= (|e_n^{(2)}|^\ast,|p_n^{(2)}|^\ast,|q_n^{(2)}|^\ast)
\end{equation*}
Then, by Lemma \ref{posmin}, the sequence $(e_n^{(3)},p_n^{(3)},q_n^{(3)})$ again forms
a minimizing sequence for $\lambda(r_1,s_1,t_1)$, and for each $n$, the real-valued functions
$e_n^{(3)}, p_n^{(3)},$ and $q_n^{(3)}$ satisfy conditions (i) and (ii) of Lemma~\ref{garlem}.
Next, it is obvious that if
$f$ and $\psi$ satisfy conditions (i) and (ii) of Lemma~\ref{garlem},
then the convolution $f\star\psi$, defined as
\begin{equation*}
\label{defconvol}
 f \star \psi (x) = \int_{-\infty}^\infty
f(x-y)\psi(y)\ dy,
\end{equation*}
also satisfies conditions (i) and (ii) of Lemma~\ref{garlem}.
Moreover, if one defines
$\psi_\epsilon = (1/\epsilon)\psi(x/\epsilon)$ for $\epsilon > 0$
and chooses $\psi$ with $\int_{-\infty}^\infty \psi(x)\ dx =
1$, then $f \star \psi_\epsilon$
converge strongly to $f$ as $\epsilon \to 0$.  Furthermore,
if $\psi\in C^\infty$ then $f \star \psi_\epsilon \in C^\infty$ as well.
Thus, by defining
\begin{equation*}
(e_n^{(4)}, p_n^{(4)},q_n^{(4)})=
(e_n^{(3)}\star\psi_{\epsilon_n},p_n^{(3)}\star\psi_{\epsilon_n},q_n^{(3)}\star\psi_{\epsilon_n}),
\end{equation*}
with $\epsilon_n$ chosen small enough for $n$ large, and setting
\begin{equation*}
f_{1, n}^{(1)}=\frac{(r_{i})^{1/2}e_{n}^{(4)}}{\|e_{n}^{(4)}\|}, \
 f_{2, n}^{(1)}=\frac{(s_{i})^{1/2}p_{n}^{(4)}}{\|p_{n}^{(4)}\|},
\text{ \
and \ }f_{3, n}^{(1)}=\frac{(t_{i})^{1/2}q_{n}^{(4)}}{\|q_{n}^{(i)}\|},
\end{equation*}
we obtain the desired minimizing sequence $(f_{1, n}^{(1)}, f_{2, n}^{(1)},f_{3, n}^{(1)})$
for $\lambda(r_1,s_1,t_1).$

\smallskip

\noindent We now proceed to prove \eqref{SUBA}. For each $n,$ let the number $x_n$ be
such that for each $1\leq j \leq 3,$ $f_{j, n}^{(1)}(x)\ \text{and} \ f_{j, n}^{(2)}(x+x_n)$
have disjoint support, and define the functions
\begin{equation*}
f_{j, n} = \left(f_{j, n}^{(1)}+\tilde f_{j, n}^{(2)}\right)^\ast, \ \textrm{where}\ \tilde f_{j, n}^{(2)}(x)=f_{j, n}^{(2)}(x+x_n), \ 1\leq j \leq 3.
\end{equation*}
It is obvious that $(f_{1,n},f_{2,n},f_{3,n})\in \Delta_{r_1+r_2,s_1+s_2,t_1+t_2}$ and hence,
\begin{equation}
\lambda(r_1+r_2,s_1+s_2,t_1 + t_2) \le \mathcal{H}(f_{1, n},f_{2, n},f_{3, n}).
\label{IleE}
\end{equation}
By making use of Lemma \ref{garlem}, it easily follows that
\begin{equation}
\begin{aligned}
\int_{-\infty}^\infty \sum_{j=1}^{3} |f_{j, nx}|^2 \ dx
& \le
\int_{-\infty}^\infty \sum_{j=1}^{3}|(f_{j, n}^{(1)}+\tilde f_{j, n}^{(2)})_x|^2 \ dx - R_n \\
 &=\int_{-\infty}^\infty \sum_{j=1}^{3}\left(|f_{j, nx}^{(1)}|^{2}+|\tilde f_{j, nx}^{(2)}|^2\right) \ dx - R_n,
\end{aligned}
\label{kinenstrictdec}
\end{equation}
where $R_n$ is given by
\begin{equation}\label{defKn}
R_n = \frac{3}{4}\sum_{j=1}^{3} \min\left\{\|f_{j, nx}^{(1)}\|^2,\|f_{j, nx}^{(2)}\|^2\right\}.
\end{equation}
Then, using the estimates \eqref{IleE}, \eqref{kinenstrictdec}, and rearrangements
properties \eqref{rearrangement1} and \eqref{rearrangement2}, we have that for all $n,$
\begin{equation}
\begin{aligned}
&\lambda(r_1+r_2, s_1+s_2, t_1+t_2) \le \mathcal{H}\left(f_{1, n},f_{2, n},f_{3, n}\right)\\
&\le \mathcal{H}\left(f_{1, n}^{(1)},f_{2, n}^{(1)},f_{3, n}^{(1)}\right)
+ \mathcal{H}\left(f_{1, n}^{(2)}, f_{2, n}^{(2)}, f_{3, n}^{(2)}\right) - R_n.
\end{aligned}
\end{equation}
Hence, by taking limit as $n\to \infty,$ we obtain
\begin{equation}
\lambda(r_1+r_2, s_1+s_2, t_1+t_2) \le \sum_{i=1}^{2} \lambda(r_i,s_i,t_i) - \liminf_{n \to \infty} R_n.
\label{subaddKn}
\end{equation}
Since $r_1+r_2>0,\ s_1+s_2>0,$ and $t_1+t_2>0,$ either both of $r_1$ and $r_2,$\ $s_1$ and $s_2,$\ $t_1$ and $t_2,$ are
positive or one of them is zero and the other is positive. To prove \eqref{SUBA}, it suffices to consider the following five cases:
\begin{enumerate}
\item[(i)] $r_1,r_2>0$ and $s_1,s_2,t_1,t_2\geq 0 ;$
\item[(ii)] $r_1=0, r_2>0, s_2>0,\ \text{and} \ t_1=0 ;$
\item[(iii)] $r_1=0, r_2>0, s_2>0, \ \text{and} \ t_1>0;$
\item[(iv)] $r_1=0, r_2>0, s_2=0, \ \text{and} \ t_1=0 ;$ and
\item[(v)] $r_1=0, r_2>0, s_2=0, \ \text{and} \ t_1>0.$
\end{enumerate}

\noindent In the case (i), i.e., when $r_1, r_2 > 0,$
Lemma
\ref{bdgnxbelow} guarantees that there exist numbers $\delta_1 >0$ and
$\delta_2 >0$ such that for all sufficiently large $n$,
$\|f_{1, nx}^{(1)}\| \ge \delta_1$ \ \text{and}\ $\|f_{1, nx}^{(2)}\| \ge
\delta_2.$
Let $\delta =
\min(\delta_1, \delta_2)> 0$. Then, \eqref{defKn} gives $R_n \ge
3\delta/4$ for all sufficiently large $n,$ and from \eqref{subaddKn},
we have
\begin{equation*}
\begin{aligned}
\lambda(r_1+r_2, s_1+s_2, t_1+t_2) &\le \lambda(r_1,s_1,t_1)+\lambda(r_2,s_2,t_2) - 3\delta/4 \\
&< \lambda(r_1,s_1,t_1)+\lambda(r_2,s_2,t_2). \label{subaddwdelta}
\end{aligned}
\end{equation*}
\noindent In the case (ii), since $r_1+s_1+t_1>0,$ so $s_1>0$ too. Then, using Lemma~\ref{bdgnxbelow} again,
there exist numbers $\delta_3, \delta_4 >0$
such that for all sufficiently large $n$,
$\|f_{2, nx}^{(1)}\| \ge \delta_3$ \ \text{and}\ $\|f_{2, nx}^{(2)}\| \ge
\delta_4.$
Let $\delta =
\min(\delta_3, \delta_4)> 0$. Then, \eqref{defKn} gives $R_n \ge
3\delta/4$ for all sufficiently large $n$ and the claim follows from \eqref{subaddKn}.

\smallskip

\noindent Next, consider the case (iii) that $r_1=0, r_2>0, s_2>0, \ \text{and} \ t_1>0.$
If $s_1>0$ or $t_2>0,$ then the proof is similar to the proof as in the case (ii) above. Thus, we may assume that
$s_1=0 \ \text{and} \ t_2=0.$ Then, we have to prove that
\begin{equation}\label{Suba3}
\lambda(r_2,s_2,t_1) < \lambda(0,0,t_1)+\lambda(r_2,s_2,0).
\end{equation}
It is well-known that (see, for example \cite{[Ca]}) the equation
\begin{equation}\label{possolnls}
-u^{\prime \prime}+\sigma_{3} u = a_{33}|u|^{2p-1}
\end{equation}
has, for any $\sigma_{3}>0,$ a unique positive solution $u_{a_{33}}$ in $H^1,$ which is explicitly given by $u_{a_{33}}(x)=e^{i\theta}\psi(x+x_0),$
where $x_0, \theta \in \mathbb{R}$ and $\psi$ is given by
\begin{equation}\label{explicitsec}
\psi(x)=\left(\frac{\sigma_3 p}{a_{33}}\right)^{1/(2p-2)} \textrm{sech}^{2/(2p-2)}\left(\frac{\sqrt{\sigma_3}(2p-2)x}{2}\right).
\end{equation}
For any $t_1>0,$ let $\psi_{t_1}$ be a solution to the problem
\begin{equation*}
\lambda(0,0,t_1) =\inf\{|f_x|_{2}^{2}-\frac{a_{33}}{p}|f|_{2p}^{2p}:\ f\in H^1 \ \text{and}\ \|f\|^{2} = t_1 \}.
\end{equation*}
Then $\psi_{t_1}$ satisfies the Lagrange multiplier equations \eqref{possolnls}, in which $\sigma_{3}$ is the Lagrange multiplier.
Therefore, $\psi_{t_1}=\psi$ up to a phase factor and a translation, where $\psi$ is
as given in \eqref{explicitsec}.
Now let $\phi_{r_2}$ and $\phi_{s_2}$ be as defined in Lemma~\ref{minforJ} so that
$\lambda(r_2,s_2,0)=F(\phi_{r_2},\phi_{s_2}).$ Then, clearly
\begin{equation*}
\int_{-\infty}^{\infty}|\phi_{r_2}|^{p} |\psi_{t_1}|^{p} >0 \ \text{and} \ \int_{-\infty}^{\infty}|\phi_{s_2}|^{p} |\psi_{t_1}|^{p} >0.
\end{equation*}
Thus, we have that
\begin{equation*}
\begin{aligned}
& \lambda(r_2,s_2,t_1)\leq \mathcal{H}(\phi_{r_2},\phi_{s_2},\psi_{t_1})=\int_{-\infty}^{\infty} \left( |(\psi_{t_1})_x|^{2}-\frac{a_{33}}{p}|\psi_{t_1}|^{2p}\right)dx \\
& + \int_{-\infty}^{\infty} \left( |(\phi_{r_2})_x|^{2}+|(\phi_{s_2})_x|^{2}-\frac{a_{11}}{p}|\phi_{r_2}|^{2p}-\frac{a_{22}}{p}|\phi_{s_2}|^{2p}
-\frac{2a_{12}}{p}|\phi_{r_2}|^{p}|\phi_{s_2}|^{p} \right)dx \\
& - \frac{2a_{13}}{p} \int_{-\infty}^{\infty}|\phi_{r_2}|^{p} |\psi_{t_1}|^{p} \ dx
- \frac{2a_{23}}{p} \int_{-\infty}^{\infty}|\phi_{s_2}|^{p} |\psi_{t_1}|^{p} \ dx,
\end{aligned}
\end{equation*}
from which it follows that
\begin{equation*}
\begin{aligned}
\lambda(r_2,s_2,t_1)& \leq \lambda(0,0,t_1) + \lambda(r_2,s_2,0)- \frac{2a_{13}}{p} \int_{-\infty}^{\infty}|\phi_{r_2}|^{p} |\psi_{t_1}|^{p} \ dx\\
&  - \frac{2a_{23}}{p} \int_{-\infty}^{\infty}|\phi_{s_2}|^{p} |\psi_{t_1}|^{p} \ dx < \lambda(0,0,t_1) + \lambda(r_2,s_2,0).\\
\end{aligned}
\end{equation*}
This proves \eqref{Suba3}. In case (iv), we have to prove that
\begin{equation}\label{Suba2}
\lambda(r_2,s_1,t_2)<\lambda(0,s_1,0) + \lambda(r_2,0,t_2),
\end{equation}
which can be proved using exactly the same argument as used in the proof of \eqref{Suba3}.
Finally, in case (v), we may assume that $t_2=0 ;$ otherwise the
claim follows from Lemma~\ref{bdgnxbelow}, \eqref{defKn}, and \eqref{subaddKn}.
Then, in case (v), we have to prove that
\begin{equation}\label{Suba4}
\lambda(r_2,s_1,t_1)<\lambda(0,s_1,t_1) + \lambda(r_2,0,0).
\end{equation}
The proof of \eqref{Suba4} is similar to the proof of \eqref{Suba3} as well.
This completes the proof of lemma. \hfill $\Box$

\section{Existence of Solitary-Wave Solutions}\label{SecVar1}

\noindent We now consider separately the three possibilities $\gamma = r+s+t,\ 0<\gamma<r+s+t,$ and $\gamma=0.$
\begin{lem}\label{GArst}
Suppose $\gamma =r+s+t.$ Then there exists a sequence of real numbers $%
\{y_{1},y_{2},y_{3},\ .\ .\ .\}$ such that

\smallskip

$1.$ for every $z<r+s+t$ there exists $\eta =\eta (z)$ such that
\begin{equation*}
\int_{y_{n}-\eta }^{y_{n}+\eta }\left(|f_{1,n}|^{2}+|f_{2,n}|^{2}+|f_{3,n}|^{2}\right) \ dx > z
\end{equation*}
for all sufficiently large $n.$

\smallskip

$2.$ the sequence $\{\mathbf{w}_{n}\}$ defined by
\begin{equation*}
w_{j,n}(x)=f_{j,n}(x+y_{n}) \ \text{for}\ x\in \mathbb{R}\  \text{and}\ j\in \{1,2,3\},
\end{equation*}
has a subsequence which converges in $Y$ norm to a function $\mathbf{\Phi} \in \mathcal{O}_{r,s,t}.$
In particular, $\mathcal{O}_{r,s,t}$ is nonempty.
\end{lem}
\noindent{\bf Proof.}
Statement 1 is just a consequence of Lions' concentration compactness lemma \cite{[L]}.
To prove statement 2, observe first that from statement 1,
there exists $\eta _{k}\in \mathbb{R}$ such
that, for every $k\in \mathbb{N},$ we have%
\begin{equation} \label{EJ10}
\int_{-\eta _{k}}^{\eta _{k}}\sum_{j=1}^{3} |w_{j, n}|
^{2} \ dx>(r+s+t)-\frac{1}{k},
\end{equation}
for all sufficiently large $n.$
As $\|w_{1,n}\|_1+\|w_{1,n}\|_1+\|w_{1,n}\|_1\leq B,$ hence
from the compact embedding of $H^{1}(\Omega)$ into $L^{2}(\Omega)$
on bounded intervals $\Omega,$ it follows that some
subsequence of $\{(w_{1,n},w_{2,n},w_{3,n})\}$ converges in $L^{2}(-\eta _{k},\eta _{k})$ norm to a limit function
$(\Phi_1,\Phi_2,\Phi_3)$ satisfying
\begin{equation*}
\int_{-\eta _{k}}^{\eta _{k}}\sum_{j=1}^{3} |\Phi_{j}|
^{2} \ dx>(r+s+t)-\frac{1}{k},
\end{equation*}
Using a Cantor diagonalization process, together with the fact that
\begin{equation}
\int_{-\infty}^{\infty}\sum_{j=1}^{3}  |w_{j, n}|
^{2} \ dx = r+s+t,\ \text{for all}\ n,
\end{equation}
we conclude that some subsequence of $\{(w_{1,n},w_{2,n},w_{3,n})\}$ converges in $L^{2}(\mathbb{R})$ norm to
a limit $(\Phi_1,\Phi_2,\Phi_3) \in L^{2}(\mathbb{R})\times L^{2}(\mathbb{R}) \times L^{2}(\mathbb{R})$ satisfying
\begin{equation*}
\int_{\infty}^{\infty}\sum_{j=1}^{3}  |\Phi_{j}|
^{2} \ dx=r+s+t.
\end{equation*}
Furthermore, by the weak
compactness of the unit sphere and the weak lower semicontinuity of the norm in Hilbert space,
$\{(w_{1,n},w_{2,n},w_{3,n})\}$ converges weakly to $(\Phi_1,\Phi_2,\Phi_3)$ in $Y,$ and that
\begin{equation*}
\|(\Phi_1,\Phi_2,\Phi_3)\|_Y \leq \liminf_{n \to \infty} \|(w_{1,n},w_{2,n},w_{3,n})\|_Y.
\end{equation*}
Next, from the Gagliardo-Nirenberg inequality, we have
\begin{equation*}
\begin{aligned}
|w_{j, n}-\Phi_j|_{2p}^{2p} & \leq C \left(\int_{-\infty}^{\infty} |w_{j,n}^{\prime}-\Phi_{j}^{\prime}|^{2} \ dx \right)^{(p-1)/2}
\left(\int_{-\infty}^{\infty} |w_{j,n}-\Phi_{j}|^{2} \ dx \right)^{(p+1)/2}\\
& \leq C \left(\int_{-\infty}^{\infty} |w_{j,n}-\Phi_{j}|^{2} \ dx \right)^{(p+1)/2},
\end{aligned}
\end{equation*}
where $C$ denotes various constants independent of $n.$ Hence $w_{j,n} \to \Phi_j$ in $L^{2p}$ norm as well.
It follows that
\begin{equation*}
\mathcal{H}(\Phi_1,\Phi_2,\Phi_3) \leq \lim_{n\to \infty} \mathcal{H}(w_{1,n},w_{2,n},w_{3,n})=\lambda(r,s,t),
\end{equation*}
whence $\mathcal{H}(\Phi_1,\Phi_2,\Phi_3)=\lambda(r,s,t)$ and $(\Phi_1,\Phi_2,\Phi_3) \in \mathcal{O}_{r,s,t}.$ Since
\begin{equation*}
|\Phi_j|_{2p} = \lim_{n \to \infty} |w_{j,n}|_{2p}, \ \|\Phi_j\|=\lim_{n \to \infty} \|w_{j,n}\|,
\end{equation*}
and
\begin{equation*}
\mathcal{H}(\Phi_1,\Phi_2,\Phi_3)=\lim_{n \to \infty} \mathcal{H}(w_{1,n},w_{2,n},w_{3,n}),
\end{equation*}
we conclude that
\begin{equation*}
\|(\Phi_1,\Phi_2,\Phi_3)\|_Y = \lim_{n \to \infty} \|(w_{1,n},w_{2,n},w_{3,n})\|_Y.
\end{equation*}
As $Y$ is a Hilbert space, an elementary
exercise in Hilbert space theory then follows that $(w_{1,n},w_{2,n},w_{3,n})$
converges to $(\Phi_1,\Phi_2,\Phi_3)$ in $Y$ norm. \hfill $\Box$

\medskip

\noindent The next result is a special case of
Lemma I.1 of \cite{[L]}. For a proof, see Lemma 2.13 of \cite{[AB11]}.
\begin{lem}  \label{Lionsvanish}
Suppose $f_n$ is a bounded sequence of functions in $H^1$
which satisfies, for some $B >0$,
\begin{equation}
\lim_{n \to \infty} \sup_{y \in \mathbb R} \int_{y-B}^{y+B} f_n^2\ dx = 0.
\label{vanishhypo}
\end{equation}
Then for every $k>2$, $|f_n|_k \to 0$ as $n \to \infty.$

\end{lem}
\noindent We can now rule out the case of vanishing:
\begin{lem}\label{GAvanish}
For any minimizing sequence $\{\mathbf{f}_{n}\}\in Y,\ \gamma > 0.$
\end{lem}
\noindent{\bf Proof.}
Suppose to contrary that $\gamma =0$. By Lemma \ref{Ibounded},
$\{|f_{1, n}|\},\ \{|f_{2, n}|\},$ and $\{|f_{3, n}|\}$
are bounded
sequences in $H^1.$ Therefore, Lemma
\ref{Lionsvanish} implies that
$|f_{j, n}|_{2p}^{2p}\ dx \to 0$ as $n\to \infty.$
For all $k, j= 1,2,3$ with $j \neq k,$ we have that
\begin{equation*}
\int_{-\infty }^{\infty }|f_{j n}|^{p}|f_{k, n}|^{p}\ dx \leq
C \left(\int_{-\infty }^{\infty }|f_{j, n}|^{2p} \ dx \right)^{1/2} \left(\int_{-\infty }^{\infty }|f_{k, n}|^{2p} \ dx\right)^{1/2},
\end{equation*}
and hence
\begin{equation*}
\lim_{n \to \infty} \int_{-\infty }^{\infty }|f_{j, n}|^{p}|f_{k, n}|^{p}\ dx =0.
\end{equation*}
It then follows that
\begin{equation}
\lambda(r, s,t)=\lim_{n\to \infty }\mathcal{H}(\mathbf{f}_{n}) \ge \liminf_{n \to
\infty}\int_{-\infty }^{\infty }\sum_{j=1}^{3} |f_{j, nx}|^2 \
dx\geq 0,
\label{istgezero}
\end{equation}
which contradicts Lemma~\ref{Ibounded}. This guarantees $\gamma > 0.$ \hfill $\Box$

\medskip

\begin{lem} \label{DIC} Suppose $r,s,t>0$ and let $\{\mathbf{f}_n\}$ be any minimizing sequence
for $\lambda(r,s,t).$ Let the number $\gamma$ be as defined in $\eqref{defgamma}.$ Then
there exist $r_1 \in [0,r], s_1 \in [0,s]$ and $t_1 \in
[0,t]$ such that
\begin{equation}
\label{SUM} \gamma = r_1+s_1+t_1
\end{equation}
and
\begin{equation} \label{REV}
\lambda(r_1, s_{1},t_{1})+\lambda(r-r_{1}, s-s_{1},t-t_{1})\leq \lambda(r,s,t).
\end{equation}
\end{lem}
\noindent{\bf Proof.}
We shall follow the arguments in \cite{[AAu]}. Let $\epsilon$ be an arbitrary positive number. By the definition of $\gamma,$ it follows
that $\gamma-\epsilon<P(\eta)\leq P(2\eta) \leq \gamma$ for $\eta$ sufficiently large.
By taking $\eta$ larger if necessary, we may also assume that $\frac{1}{\eta}<\epsilon.$
From the definition of $P,$ we can choose $N$ so large that, for every $n\geq N,$
\begin{equation*}
\gamma-\epsilon < P_{n}(\eta) \leq P_{n}(2\eta) \leq \gamma +\epsilon.
\end{equation*}
Hence, for each $n\geq N,$ we can find $y_n$ such that
\begin{equation}\label{sub113}
\int_{y_{n-\eta}}^{y_{n}+\eta}\sum_{j=1}^{3}|f_{j, n}|
^{2} \ dx > \gamma -\epsilon \ \text{and} \ %
\int_{y_{n-2\eta}}^{y_{n}+2\eta}\sum_{j=1}^{3}|f_{j, n}|
^{2} \ dx < \gamma +\epsilon.
\end{equation}
Now choose $\rho \in C_{0}^{\infty }[-2,2]$ such that $\rho \equiv 1$ on $[-1,1],$ and
let $\sigma \in C^{\infty}(\mathbb{R})$ be such that
$\rho^2 + \sigma^2 \equiv 1$ on $\mathbb R$.
Set, for $\eta > 0,$
$\rho_\eta(x)=\rho(x/\eta)$ \ \text{and} \ $\sigma_\eta(x)=\sigma(x/\eta).$
and define the functions
\begin{equation*}
\mathbf{f}_{n}^{(1)}(x)=\rho_\eta(x-y_n)\mathbf{f}(x) \ \text{and} \ \mathbf{f}_{n}^{(2)}(x)=\sigma_\eta(x-y_n)\mathbf{f}(x).
\end{equation*}
Then, for each $j=1,2,3,$ and $k=1,2,$ the sequences $\{f_{j, n}^{(k)}\}$ are bounded in $L^{2}.$ Thus,
by passing to subsequences if necessary, we may assume that there exist
$r_1 \in [0,r], s_1 \in [0,s]$ and $t_1 \in
[0,t]$ such that
\begin{equation}\label{Aux3}
\|f_{1, n}^{(1)}\|^{2} \to r_1, \ \|f_{2, n}^{(1)}\|^{2} \to s_1, \ \text{and} \
\|f_{3, n}^{(1)}\|^{2} \to t_1,
\end{equation}
whence it follows also that
\begin{equation}\label{Aux4}
\|f_{1, n}^{(2)}\|^{2} \to r-r_1, \ \|f_{2, n}^{(2)}\|^{2} \to s-s_1, \ \text{and} \
\|f_{3, n}^{(2)}\|^{2}\to t-t_1.
\end{equation}
Now
\begin{equation*}
r_1+s_{1}+t_{1}=\lim_{n\to \infty }\int_{-\infty
}^{\infty } \sum_{j=1}^{3}|f_{j, n}^{(1)}|^{2}\ dx=\lim_{n\to \infty
}\int_{-\infty }^{\infty }\rho _{\eta }^{2}\sum_{j=1}^{3}|
f_{j, n}|^{2} \ dx.
\end{equation*}
From \eqref{sub113}, it follows that, for every $n\in N,$
\begin{equation*}
\gamma -\epsilon <\int_{-\infty }^{\infty }\rho _{\eta }^{2}\sum_{j=1}^{3}|
f_{j, n}|^{2} \ dx <\gamma +\epsilon .
\end{equation*}
Hence $|(r_1+s_{1}+t_{1})-\gamma| <\epsilon .$
Next, we claim that for all $n,$
\begin{equation} \label{e12ineq}
\mathcal{H}(\mathbf{f}_{n}^{(1)})+\mathcal{H}(\mathbf{f}_{n}^{(2)}) \le \mathcal{H}(\mathbf{f}_n) + C\epsilon
\end{equation}
\noindent To see \eqref{e12ineq}, we write
\begin{equation*}
\begin{aligned}
\mathcal{H}(\mathbf{f}_{n}^{(1)}) &= \int_{-\infty }^{\infty }\rho _{\eta }^{2}\left(
\sum_{j=1}^{3}|f_{j,nx}|^{2} - \frac{1}{p} \sum_{k,j=1}^{3} a_{kj} |f_{k, n}|^{p}|f_{j, n}|^{p} \right) \ dx\\
&\ \  +\int_{-\infty }^{\infty }\left(\left( \rho _{\eta }^{\prime
}\right) ^{2}\sum_{j=1}^{3} |f_{j, n}|^{2} +2\rho _{\eta }^{\prime
}\rho _{\eta } \sum_{j=1}^{3} f_{j, n} f_{j, nx} \right) \ dx \\
&\ \ +\frac{1}{p} \int_{-\infty }^{\infty } (\rho_{\eta}^{2}-\rho_{\eta}^{2p}) \sum_{k,j=1}^{3} a_{kj} |f_{k, n}|^{p}|f_{j, n}|^{p} dx,
\end{aligned}
\end{equation*}
where, for ease of notation, we have written the functions $\rho_\eta (x-y_n)$ simply as $\rho_{\eta}.$
Similar estimate holds for $\mathcal{H}(\mathbf{f}_{n}^{(2)}).$
Since $\rho_{\eta}^2 + \sigma_{\eta}^2 \equiv 1, \ |\rho_{\eta}^{\prime}|_\infty=|\rho^{\prime}|_\infty / \eta,$ and
$|\sigma_{\eta}^{\prime}|_\infty=|\sigma^{\prime}|_\infty / \eta,$ an application of H\"{o}lder's inequality yields
\begin{equation*}
\begin{aligned}
\mathcal{H}(\mathbf{f}_{n}^{(1)})& +\mathcal{H}(\mathbf{f}_{n}^{(2)}) =\mathcal{H}(\mathbf{f}_{n})+O(1/\eta) \\
\ &+\frac{1}{p} \int_{-\infty }^{\infty } \left[(\rho_{\eta}^{2}-\rho_{\eta}^{2p})
 +(\sigma_{\eta}^{2}-\sigma_{\eta}^{2p}) \right] \sum_{k,j=1}^{3} a_{kj} |f_{k, n}|^{p}|f_{j, n}|^{p} dx,
\end{aligned}
\end{equation*}
where $O(1/\eta)$ denotes a term bounded in absolute value by $C/\eta$ with $C$ independent of
$\eta$ and $n.$ Using \eqref{sub113}, one can see that
\begin{equation*}
\begin{aligned}
\left\vert \int_{-\infty }^{\infty }\left[ (\rho _{\eta }^{2}-\rho
_{\eta }^{2p})+(\sigma _{\eta }^{2}-\sigma _{\eta }^{2p})\right] |
f_{k, n} f_{j, n}|^{p}\ dx \right\vert
 & \leq |f_{k, n}|_{\infty }^{p}\int_{\eta \leq
|x-y_{n}| \leq 2\eta }2|f_{j, n}|^{p} \ dx \\
& \leq C\epsilon,
\end{aligned}
\end{equation*}
where again $C$ denotes various constants independent of $\eta$ and $n.$
Then, \eqref{e12ineq} follows by choosing $\eta$ large enough so that $|O(1/\eta)|\leq \epsilon.$

\smallskip

\noindent To prove \eqref{REV}, notice that for any given value of $\epsilon,$ each of the terms in \eqref{e12ineq}
is bounded independently of $n,$ so by passing to a subsequence if necessary, we may assume that
\begin{equation}\label{Aux5}
\mathcal{H}(\mathbf{f}_{n}^{(1)})\to H_1 \  \text{and} \ \mathcal{H}(\mathbf{f}_{n}^{(2)}) \to H_2.
\end{equation}
Then,
$H_1 +H_2 \leq \lambda(r,s,t)+C\epsilon.$
Since $\epsilon$ can be taken arbitrarily small and $\eta$ arbitrarily large, combining the results of the preceding paragraphs,
we can find sequences $\{\mathbf{f}_{n}^{(1, k)}\}$ and $\{\mathbf{f}_{n}^{(2, k)}\},$ for each $k\in \mathbb{N},$ such that
\begin{equation*}
\begin{aligned}
& \|f_{1,n}^{(1, k)}\|^{2}\to r_1(k), \ \|f_{2,n}^{(1, k)}\|^{2}\to s_1(k), \ \|f_{3,n}^{(1, k)}\|^{2}\to t_1(k), \\
& \|f_{1,n}^{(2, k)}\|^{2}\to r-r_1(k), \ \|f_{2,n}^{(2, k)}\|^{2}\to s-s_1(k), \ \|f_{3,n}^{(2, k)}\|^{2}\to t-t_1(k), \\
& \mathcal{H}(\mathbf{f}_{n}^{(1, k)})\to H_1(k), \  \text{and} \ \mathcal{H}(\mathbf{f}_{n}^{(2, k)}) \to H_2(k),
\end{aligned}
\end{equation*}
where $r_1(k) \in [0,r], \ s_1(k) \in [0,s], \ t_1(k) \in [0,t],$
\begin{equation}\label{Aux1}
|r_1(k)+s_1(k)+t_1(k)-\gamma| \leq \epsilon,
\end{equation}
and
\begin{equation}\label{Aux2}
H_1(k) +H_2(k) \leq \lambda(r,s,t)+\frac{1}{k}.
\end{equation}
By passing to subsequences, we may assume that
\begin{equation*}
\begin{aligned}
& r_1(k) \to r_1\in [0,r], \ s_1(k) \to s_1\in [0,s], \ t_1(k) \to t_1 \in [0,t], \\
& H_1(k) \to H_1, \ \text{and} \ H_2(k)\to H_2.
\end{aligned}
\end{equation*}
Also, by redefining $\{\mathbf{f}_{n}^{(1)}\}$ and $\{\mathbf{f}_{n}^{(2)}\}$ as the diagonal subsequences
\begin{equation*}
\mathbf{f}_{n}^{(1)}=\mathbf{f}_{n}^{(1, n)} \ \text{and} \ \mathbf{f}_{n}^{(2)}=\mathbf{f}_{n}^{(2, n)},
\end{equation*}
we may assume that \eqref{Aux3}, \eqref{Aux4}, and \eqref{Aux5} hold.

\smallskip

\noindent By letting $k \to \infty$ in \eqref{Aux1} yields \eqref{SUM}. The claim \eqref{REV}
follows from \eqref{Aux2} provided we can show that
\begin{equation}\label{Aux6}
H_1 \geq \lambda(r_1, s_1, t_1),
\end{equation}
and
\begin{equation}\label{Aux7}
H_2 \geq \lambda(r-r_1, s-s_1, t-t_1).
\end{equation}
To prove \eqref{Aux6}, consider first the case that $r_1, s_{1},$ and $t_{1}$ are all positive.
Then, for $n$ sufficiently large,
 $\|f_{j, n}^{(1)}\|$ are all positive for each $j=1,2,3,$ so we may define
\begin{equation*}
\alpha _{n}=\frac{\sqrt{r_{1}}}{\| f_{1, n}^{(1)}\| },\ \beta
_{n}=\frac{\sqrt{s_{1}}}{\| f_{2, n}^{(1)}\| },\ \text{and} \ \gamma _{n}=%
\frac{\sqrt{t_1}}{\| f_{3, n}^{(1)}\| },
\end{equation*}
which gives
$(\alpha _{n}f_{1, n}^{(1)}, \beta _{n}f_{2, n}^{(1)}, \gamma_{n}f_{3, n}^{(1)}) \in \Delta_{r_1,s_1,t_1}.$ Consequently, we have
\begin{equation*}
\mathcal{H}(\alpha _{n}f_{1, n}^{(1)}, \beta _{n}f_{2, n}^{(1)}, \gamma_{n}f_{3, n}^{(1)}) \geq \lambda(r_1, s_1, t_1).
\end{equation*}
As all scaling factors tend to $1$ as $n\to \infty ,$ it follows that
\begin{equation*}
\mathcal{H}(\alpha _{n}f_{1, n}^{(1)}, \beta _{n}f_{2, n}^{(1)}, \gamma_{n}f_{3, n}^{(1)}) \to H_1,
\end{equation*}
and hence \eqref{Aux6} follows. Next, we prove \eqref{Aux6} if exactly one of $r_1, \ s_1,$ or $t_1$ is zero.
Consider the case that $r_1=0, s_1>0,$ and $t_1>0 .$ Then, using the Gagliardo-Nirenberg inequality, we have that
\begin{equation*}
\int_{-\infty}^{\infty} |f_{1, n}|^{p} |f_{j, n}|^{p} \ dx \to 0 \ \ \text{for all} \ 1\leq j \leq 3,
\end{equation*}
and hence, we deduce that
\begin{equation*}
\begin{aligned}
H_{1}& =\lim_{n\to \infty}\mathcal{H}(\mathbf{f}_{n}^{(1)})=\lim_{n\to
\infty}\int_{-\infty }^{\infty }\left(|f_{1,nx}^{(1)}|^{2}+|f_{2,nx}^{(1)}|^{2}+|f_{3,nx}^{(1)}|^{2}\right. \\
& -\frac{a_{22}}{p}|f_{2, n}^{(1)}|^{2p} - \frac{a_{33}}{p}
|f_{3, n}^{(1)}|^{2p}-\left. \frac{2a_{23}}{p}|f_{2, n}^{(1)}|^{p}|f_{3, n}^{(1)}|^{p} \right) \ dx \\
& \geq \liminf_{n\to \infty}\int_{-\infty }^{\infty }\left( \sum_{j=2}^{3}|f_{j,nx}^{(1)}|^{2} - \frac{1}{p}
\sum_{k,j=2}^{3}a_{kj}|f_{k, n}^{(1)}|^{p}|f_{j, n}^{(1)}|^{p}\right) dx \geq \lambda(0, s_1, t_1).
\end{aligned}
\end{equation*}
Finally, if $r_1=0, s_1=0,$ and $t_1>0,$ then we have
\begin{equation*}
\begin{aligned}
H_{1} & =\lim_{n\to \infty}\mathcal{H}(\mathbf{f}_{n}^{(1)})=\lim_{n\to
\infty}\int_{-\infty }^{\infty }\left( \sum_{j=1}^{3}|f_{j,nx}^{(1)}|^{2} -\frac{a_{33}}{p}|f_{3, n}^{(1)}|^{2p}\right) \ dx \\
& \geq \liminf_{n\to \infty}\int_{-\infty }^{\infty } \left( |f_{3,nx}^{(1)}|^{2} -\frac{a_{33}}{p}|f_{3, n}^{(1)}|^{2p} \right)
\geq \lambda(0,0,t_1).
\end{aligned}
\end{equation*}
This completes the proof of \eqref{Aux6}. The proof of \eqref{Aux7} is similar with $r-r_1, s-s_1,$ and $t-t_1$ playing the
roles of $r_1, s_1,$ and $t_1,$ respectively. \hfill $\Box$

\medskip

\noindent The following lemma rules out the possibility of dichotomy of minimizing sequences:
\begin{lem}\label{GAdicho}
For every minimizing sequence, $\gamma \not\in (0, r+s+t).$
\end{lem}
\noindent{\bf Proof.}
We proceed by contradiction. Suppose that $0<\gamma < r+s + t$.   Let $r_1, s_1,$ and
$t_1$ be as in Lemma~\ref{DIC}, and define $r_2=r-r_1, s_2=s-s_1,$ and
$t_2=t-t_1$. It then follows that $r_2+s_{2}+t_{2}=(r+s+t)-\gamma >0$, and also
$r_1+s_{1}+t_{1}=\gamma >0$.
Furthermore, $r_1+r_2=r>0,\ s_{1}+s_{2}=s>0,$ and $t_{1}+t_{2}=t>0$. Therefore, as a consequence of
Lemma~\ref{subadd}, \eqref{SUBA} holds. But this
contradicts the fact \eqref{REV} and thus, lemma follows. \hfill $\Box$

\medskip

\noindent The next theorem guarantees the existence of a minimizing pair for
\eqref{Ist} and hence, the existence of three-parameter family of solitary waves for the 3-coupled NLS system \eqref{CNLS}
provided that $a_{kj}>0$ for all $k,j\in \{1,2,3\}$ and all $2 \leq p < 3.$

\begin{thm}\label{Main1}
The set $\mathcal{O}_{r,s,t}$ is not empty. Moreover, if $\{\mathbf{f}_{n}\}$ is any minimizing sequence
for $\lambda(r,s,t)$, then

\smallskip

1. There exists a sequence $\{y_{k}\} \subset \mathbb{R}$
and an element $\mathbf{\Phi} \in \mathcal{O}_{r,s,t}$ such that
$\{\mathbf{f}_{n}(\cdot +y_{n})\}$ has a subsequence converging strongly in $Y$ to $\mathbf{\Phi}.$

\smallskip

2. Each function
$\mathbf{\Phi}\in\mathcal{O}_{r,s,t}$ is a solution of the system \eqref{ODE}
for some $\omega_{1}, \omega_{2}, \omega_{3} >0$, and therefore when inserted into \eqref{SO} yields
a three parameter family solitary-wave
solution of the NLS system \eqref{CNLS}.

\smallskip

3. The following holds:
\begin{equation*}
\lim_{n\to \infty} \inf_{y \in \mathbf{R}} \inf_{\mathbf{\Phi} \in \mathcal{O}_{r,s,t}} \|\mathbf{f}_{n}(\cdot +y) - \mathbf{\Phi}\|_{Y} =0.
\end{equation*}

\smallskip

4. The following holds:
\begin{equation*}
\lim_{n\to \infty} \inf_{\mathbf{\Phi} \in \mathcal{O}_{r,s,t}} \|\mathbf{f}_{n} - \mathbf{\Phi}\|_{Y} =0.
\end{equation*}
\end{thm}
\noindent{\bf Proof.}
From Lemmas~\ref{GAvanish} and \ref{GAdicho}, it follows that $\gamma=r+s+t.$ Then,
by Lemma~\ref{GArst}, the set $\mathcal{O}_{r,s,t}$ is not empty and statement 1 follows.

\smallskip

\noindent To prove statement 2, since $\mathbf{\Phi}$ is a minimizing function for $\lambda(r,s,t)$ and so, using the
Lagrange multiplier principle,
there exist real numbers $\omega_{1},\ \omega_2,$ and $\omega_{3}$ such that
\begin{equation}
\delta \mathcal{H}(\Phi_1,\Phi_2, \Phi_3) + \omega_{1} \delta \mathcal{Q}(\Phi_1)
+\omega_{2} \delta \mathcal{Q}(\Phi_2)+\omega_3 \delta \mathcal{Q}(\Phi_3)=0.
\end{equation}
One can now see by computing the
Fr\'echet derivatives that the equations
\begin{equation}\label{ODE1}
\left\{
\begin{aligned}
-\Phi_{1, xx}+\omega_{1}\Phi_1 &=a_{11}|\Phi_1|^{2p-2}\Phi_1+\left(a_{12}|\Phi_2|^{p}+a_{13}|\Phi_3|^{p}\right)|\Phi_1|^{p-2}\Phi_1, \\
-\Phi_{2, xx}+\omega_{2}\Phi_2 &=a_{22}|\Phi_2|^{2p-2}\Phi_2+\left(a_{12}|\Phi_1|^{p}+a_{23}|\Phi_3|^{p}\right)|\Phi_2|^{p-2}\Phi_2, \\
-\Phi_{3, xx}+\omega_{3}\Phi_3 &=a_{33}|\Phi_3|^{2p-2}\Phi_3+\left(a_{13}|\Phi_1|^{p}+a_{23}|\Phi_2|^{p}\right)|\Phi_3|^{p-2}\Phi_3,
\end{aligned}
\right.
\end{equation}
hold (in distributional sense). A straightforward bootstrapping argument (for example, \ Lemma 1.3 of Tao's book \cite{[TTao]}) shows that these
distributional
solutions are in fact classical solutions.

\smallskip

\noindent Multiplying the first equation in \eqref{ODE1} by $\bar{\Phi}_1, $ the second
equation by $\bar{\Phi}_2 ,$ and the third equation by $\bar{\Phi}_3 ,$ and integrating over $\mathbb{R},$ we obtain that
\begin{equation}\label{ODE2}
\int_{-\infty }^{\infty }\left(
|%
\Phi_{j}^{\prime}|^{2}-|\Phi_j|^p \sum_{k=1}^{3}a_{jk}
|\Phi_k|^{p} \right)\ dx = -\omega_j \int_{-\infty }^{\infty } |\Phi_j|^{2}\ dx, \ j=1, 2, 3.
\end{equation}
By Lemma~\ref{negative}, with $\mathbf{f}_{n}=\mathbf{\Phi},$ one has that
\begin{equation*}
\int_{-\infty }^{\infty }\left(
|\Phi_{j}^{\prime}|^{2}-|\Phi_j|^p \sum_{k=1}^{3}a_{jk}
|\Phi_k|^{p} \right)\ dx <0, \ j=1,2,3,
\end{equation*}
and hence, $\omega_{1}, \omega_{2}, \omega_3>0.$ This
then completes the proof of statement 2.

\smallskip

\noindent To prove statement 3, suppose that it is false. Then, there there exists a subsequence $\{\mathbf{f}_{n_k}\}$ of $\{\mathbf{f}_{n}\}$
and a number $\varepsilon>0$ such that
\begin{equation*}
\lim_{n\to \infty} \inf_{y\in \mathbb{R}} \inf_{\mathbf{\Phi} \in \mathcal{O}_{r,s,t}}
\|\mathbf{f}_{n}(\cdot +y) - \mathbf{\Phi}\|_{Y} \geq \varepsilon
\end{equation*}
for all $k \in \mathbb{N}.$
As $\{\mathbf{f}_{n_k}\}$ itself a minimizing sequence for $\lambda(r,s,t),$ it
follows from statement 1 that there exists
a sequence of real numbers $\{y_k\}$ and an
element $(\Phi_{1,0},\Phi_{2,0},\Phi_{3,0})$ of $\mathcal{O}_{r,s,t}$ such that
\begin{equation*}
\liminf_{k\to \infty} \|\mathbf{f}_{n_k}(\cdot +y_k) - (\Phi_{1,0},\Phi_{2,0},\Phi_{3,0})\|_{Y} =0.
\end{equation*}
This contradiction proves statement 3.

\smallskip

\noindent Finally, since $\mathcal{H}$ and $\mathcal{Q}$ are invariant under translations, then $\mathcal{O}_{r,s,t}$
clearly contains any translate of $\mathbf{\Phi}$ if it contains $\mathbf{\Phi},$
and hence, statement 4 follows from statement 3. \hfill $\Box$

\medskip

\noindent The next theorem addresses the question about the characterization of the
set $\mathcal{O}_{r,s,t}.$

\begin{thm}\label{Main2}
For every $\mathbf{\Phi}$ in $\mathcal{O}_{r,s,t}$, there exist numbers
$\theta_{j} \in \mathbb{R}$ and real functions $\phi_{j}$ such that $\phi_{j}(x) >0,$
for all $x \in \mathbb R$, and
\begin{equation*}
\Phi_j(x)= e^{i\theta_{j}}\phi_{j}(x), \ j=1,2,3.
\end{equation*}
Furthermore,
$\Phi_{1}, \Phi_2, \Phi_3$ are infinitely differentiable on $\mathbb R$.
\end{thm}
\noindent{\bf Proof.}
We write the complex-valued functions $\Phi_{j}$ as
\begin{equation}\label{char1}
\Phi_{j} (x)=e^{i\theta _{j}(x)} \phi_{j}(x), \ j=1,2,3 ,
\end{equation}
where $\theta _{j}:\mathbb{R}\to \mathbb{R}$ and $\phi_{j}(x)=|
\Phi_{j}(x)|, \ j=1,2,3. $ Notice that $(\phi_1 ,\phi_2,\phi_3 )$ is also
in $\mathcal O_{r,s,t}$, as follows from Lemma~\ref{posmin}. Therefore, $(\phi_1 ,\phi_2,\phi_3 )$
satisfies the Lagrange multiplier equations
\begin{equation}\label{ODE4}
-\phi_{j, xx}+\omega_{j}\phi_j =|\phi_j|^{p-2}\phi_j\sum_{k=1}^{3}a_{jk}|\phi_k |^{p},\ j=1,2,3.
\end{equation}
(That the Lagrange multipliers stay same follows from the fact that they are
determined by \eqref{ODE2}, and this equation remains unchanged when one replaces $(\Phi_1 ,\Phi_2, \Phi_3
)$ by $(\phi_1 ,\phi_2,\phi_3 )$.)
Using \eqref{char1}, we now compute
\begin{equation}\label{pos1}
\Phi_{1}^{\prime \prime}=e^{i\theta _{1}(x)}\left( \omega_1 \phi_1
-|\phi_1|^{p-2}\phi_1\sum_{k=1}^{3}a_{1k}|\phi_k|^{p} -Z(x)
 \right),
\end{equation}
where
\begin{equation*}
Z(x)=(\theta _{1}^{\prime }(x))^{2}
\phi_1(x) - 2i\theta _{1}^{\prime }(x)\phi_{1}^{\prime }(x) - i\theta
_{1}^{\prime \prime }(x) \phi_1(x).
\end{equation*}
On the other hand, since $(\Phi_1, \Phi_2, \Phi_3)$ satisfies the same equations \eqref{ODE4} as $(\phi_1, \phi_2, \phi_3)$, it follows that
\begin{equation}\label{pos2}
\Phi_{1} ^{\prime \prime }=e^{i\theta _{1}(x)}\left(\omega_1 \phi_1-|\phi_1|^{p-2}\phi_1 \sum_{k=1}^{3}a_{1k}|\phi_k|^{p} \right).
\end{equation}

\noindent From \eqref{pos1} and \eqref{pos2} , we obtain that $Z(x)=0,$ and by
equating the real part of this equation, we conclude that $\theta
_{1}^{\prime }(x)=0,$ and hence $\theta _{1}(x)$ is constant. Similarly, $%
\theta _{2}(x)$ and $\theta _{3}(x)$ are constants.

\smallskip

\noindent Next, a straightforward calculation using Fourier transform shows that for each $j=1,2,3,$
the operator $-\partial_{x}^{2}+\omega_j$ appearing in \eqref{ODE4} is invertible on $H^1,$ with
inverse given by convolution with
the function
\begin{equation*}
E_{\omega_j}(x)=\frac{1}{2\sqrt{\omega_j}}e^{-\sqrt{\omega_j} |x|}.
\end{equation*}
The Lagrange multiplier equations associated with $(\phi_1,\phi_2,\phi_3)$ can then be rewritten in the form
\begin{equation*}
\begin{aligned}
& \phi_1=E_{\omega_1}
\star\left(a_{11}|\phi_1|^{2p-2}\phi_1+a_{12}|\phi_2|^p |\phi_1|^{p-2} \phi_1 +a_{13}|\phi_3|^p |\phi_1|^{p-2} \phi_1 \right),\\
& \phi_2=E_{\omega_2}
\star\left(a_{22}|\phi_2|^{2p-2}\phi_2+a_{12}|\phi_1|^p |\phi_2|^{p-2} \phi_2 +a_{23}|\phi_3|^p |\phi_2|^{p-2} \phi_2 \right),\\
& \phi_3=E_{\omega_3}
\star\left(a_{33}|\phi_3|^{2p-2}\phi_3+a_{13}|\phi_1|^2 |\phi_3|^{p-2} \phi_3 +a_{23}|\phi_2|^p |\phi_3|^{p-2} \phi_3 \right).
\end{aligned}
\end{equation*}
Since the convolutions of the positive kernel $E_{\omega_j}$ with functions which are not identically zero and non-negative 
everywhere on $\mathbb{R}$
produce everywhere positive functions on $\mathbb{R}$, we conclude
that $\phi_{j}(x)>0$ on $\mathbb R.$ \hfill $\Box$

\section{Stability of Solitary Waves}

Our stability result reads as follows.

\begin{thm}\label{Main3}
For every $\epsilon >0$, there
exists $\delta>0$ such that if
\begin{equation*}
\inf_{\mathbf{\Phi}\in \mathcal{O}_{r,s,t}}\|(f_0, g_0, h_0)-\mathbf{\Phi}\|_{Y}<\delta ,
\end{equation*}
then the solution $\mathbf{u}(x,t)$
of \eqref{CNLS} with
$\mathbf{u}(x,0) =(f_0(x),g_0(x),h_0(x))$ satisfies
\begin{equation*}
\sup_{t \in \mathbb{R}}\inf_{\mathbf{\Phi}\in \mathcal{O}_{r,s,t}}\|\mathbf{u}(\cdot,t)-\mathbf{\Phi}\|_{Y}<\epsilon
\end{equation*}
 \end{thm}
\noindent{\bf Proof.} The proof follows a standard argument.
Suppose that the set $\mathcal{O}_{r,s,t}$ is not stable. Then there exist a
number $\epsilon >0,$ a sequence of times ${t_{n}},$ and a sequence
$\{\mathbf{u}_n(x,0)\}$ in $Y$ such that for all $n,$
\begin{equation}
\label{idconvtoF} \inf\{\|\mathbf{u}_n(x,0)-\mathbf{\Phi} \|_Y : \mathbf{\Phi} \in
\mathcal{O}_{r,s,t}\}<\frac{1}{n};
\end{equation}
and
\begin{equation}\label{contra}
\inf\{\|\mathbf{u}_n(\cdot,t_n)-\mathbf{\Phi} \|_Y : \mathbf{\Phi} \in
\mathcal{O}_{r,s,t}\} \geq \epsilon,
\end{equation}
for all $n,$ where $\mathbf{u}_n(x,t)$ solves \eqref{CNLS} with initial data $\mathbf{u}_n(x,0).$
Since $\mathbf{u}_n(x,0)$ converges to an element in $\mathcal{O}_{r,s,t}$ in $Y$ norm, and since
for $\mathbf{\Phi} \in \mathcal{O}_{r,s,t},$ we have
$\mathcal{Q}(\Phi_1)=r, \ \mathcal{Q}(\Phi_2)=s, \ \mathcal{Q}(\Phi_3)=t, \ \text{and} \ \mathcal{H}(\mathbf{\Phi})=\lambda(r,s,t),$
we therefore have
\begin{equation*}
\mathcal{Q}(u_{1, n}(x,0))\to r, \ \mathcal{Q}(u_{2, n}(x,0))\to s, \ \mathcal{Q}(u_{3, n}(x,0))\to t,
\end{equation*}
and $\mathcal{H}(\mathbf{u}_{ n}(x,0)) \to \lambda(r,s,t).$
Let us denote $u_{1, n}(\cdot,t_n)$ by $U_{1, n}, u_{2, n}(\cdot,t_n)$ by $U_{2, n},$
and $u_{3,n}(\cdot,t_n)$ by $U_{3,n}$.
We now choose $\{\alpha_n\}, \{\beta_n\},$ and $\{\gamma_n\}$ such that
\begin{equation*}
\mathcal{Q}(\alpha_n u_{1, n}(x,0))= r, \ \mathcal{Q}(\beta_n u_{2, n}(x,0))= s, \ \mathcal{Q}(\gamma_n u_{3, n}(x,0))= t,
\end{equation*}
for all $n.$ Thus, $\alpha_n \to 1, \beta_n \to 1,$ and $ \gamma_n \to 1.$
Hence the sequences $f_{1,n}=\alpha_n U_{1, n}, f_{2,n}=\beta_n U_{2, n},$ and
$f_{3,n}=\gamma_n U_{3, n}$ satisfies
$\mathcal{Q}(f_{1,n})=r, \mathcal{Q}(f_{2,n})=s, \mathcal{Q}(f_{3,n})=t,$ and
\begin{equation*}
\lim_{n\to \infty }\mathcal{H}(\mathbf{f}_{n})= \lim_{n\to \infty }\mathcal{H}(\mathbf{u}_{n}(\cdot, t_n))
 = \lim_{n\to \infty }\mathcal{H}(\mathbf{u}_{n}(x, 0))=\lambda(r,s,t).
\end{equation*}
Therefore, $\{\mathbf{f}_n\}$ is a minimizing sequence for $\lambda(r,s,t).$ From Theorem~\ref{Main1}, it follows that
for all $n$ sufficiently large, there exists $\mathbf{\Phi}_n \in \mathcal{O}_{r,s,t}$ such that
$\|\mathbf{f}_{n}-\mathbf{\Phi}_n\|_{Y} < \epsilon /2.$ But then we have
\begin{equation*}
\begin{aligned}
\epsilon & \leq \|\mathbf{u}_{n}(\cdot, t_n)-\mathbf{\Phi}_n\|_{Y} \leq \|\mathbf{u}_{n}(\cdot, t_n)-\mathbf{f}_n \|_{Y}
+ \|\mathbf{f}_n -\mathbf{\Phi}_n\|_{Y} \\
& \leq |1-\alpha_n|\cdot \|U_{1, n}\|_1 +|1-\beta_n|\cdot \|U_{2, n}\|_1
+ |1-\gamma_n|\cdot \|U_{3, n}\|_1 + \frac{\epsilon}{2}
\end{aligned}
\end{equation*}
and by taking $n \to \infty,$ we obtain that $\epsilon \leq \epsilon /2,$ a contradiction,
 and we conclude that $\mathcal
O_{r,s,t}$ must in fact be stable. \hfill $\Box$

\medskip

\noindent \textsc{Trocaire College, 360 Choate Ave, Buffalo, NY 14220}

\smallskip

\noindent \textit{E-mail address:} \texttt{bhattarais@trocaire.edu, sntbhattarai@gmail.com}

\end{document}